\documentclass[11pt]{amsart}

\usepackage[a4paper,margin=1.05in]{geometry}
\usepackage{latexsym,rawfonts}
\usepackage{amsfonts,amssymb}
\usepackage{amsmath,amsthm}

\usepackage[plainpages=false]{hyperref}
\usepackage{graphicx}
\usepackage{color}

\newtheorem{theorem}{Theorem}[section]

\newtheorem{lemma}[theorem]{Lemma}
\newtheorem{proposition}[theorem]{Proposition}

\theoremstyle{definition}

\numberwithin{equation}{section}

\DeclareMathOperator{\TR}{tr}

\DeclareMathOperator{\SUPP}{supp}
\DeclareMathOperator{\DIAM}{diam}
\DeclareMathOperator{\DIAG}{diag}

\newcommand{\R}{\mathbb{R}}

\begin{document}

\title{Global regularities for the Dirichlet problem of Monge-Amp\`ere type equations}

\author{Feida Jiang}
\address{School of Mathematics and Shing-Tung Yau Center, Southeast University, Nanjing 211189, China}
\email{jiangfeida@seu.edu.cn}

\author{Jiakun Liu}
\address{School of Mathematics and Statistics, The University of Sydney, Camperdown, NSW 2006,
Australia}
\email{jiakun.liu@sydney.edu.au}

\author{Jian Lu}
\address{School of Mathematics and Statistics, Nanjing University of Science and Technology, 
  Nanjing 210094, China}
\email{lj-tshu04@163.com}

\author{Chuanlong Sun}
\address{School of Mathematics and Shing-Tung Yau Center, Southeast University, Nanjing 211189, China}
\email{cl-sun@seu.edu.cn}

\date{\today}
\expandafter\providecommand\csname subjclassname@2020\endcsname
  {\textup{2020} Mathematics Subject Classification}
\begin{abstract}
We establish global a priori $W^{2,p}$ and $C^{2,\alpha}$ estimates for elliptic solutions of the Dirichlet problem for cost-generated Monge-Amp\`ere type equations under the weak Ma-Trudinger-Wang condition A3w, assuming continuous and H\"older continuous right-hand sides, respectively.
\end{abstract}
\maketitle

 \baselineskip15pt
 \parskip4pt

\section{Introduction}

In this paper we consider the Dirichlet problem
	\begin{equation}\label{eq:problem}
	\begin{split}
		 \det \left[ D^2u-A(x,Du) \right] &= f(x) \quad\hbox{in }\Omega,\\
		 		 u &= \varphi \qquad\hbox{ on }\partial\Omega,
	\end{split}
	\end{equation}
where the matrix $A$ is generated by a cost function through
	\begin{equation}\label{eq:cost}
 		A(x,p)=c_{xx}(x,Y(x,p)),\qquad c_x(x,Y(x,p))=p.
	\end{equation}
The equation in \eqref{eq:problem} arises in optimal transportation, where it is coupled with a second boundary condition prescribing the image of the transport map $T_u(x)=Y(x,Du(x))$; see \cite{DF14,MTW05,TrudingerWang2009,Villa1,Villa2}.
Here, we study the Dirichlet problem building on our previous work \cite{HuangJiangLiu2015}. 
We adopt the sign convention that a $c$-support has the form $c(x,y)+a$ and
lies below $u$ throughout $\Omega$, so that $u$ is $c$-convex. 

Boundary regularity for the classical Monge-Amp\`ere equation, corresponding to the cost $c(x,y)=x\cdot y$, provides a geometric model for this problem, see \cite{TrudingerWang2008,Savin2013Boundary,Savin2013W2p}.
For smooth $f$, Jiang, Trudinger and Yang
\cite{JiangTrudingerYang2014} established global second derivative
estimates under A3w by assuming the existence of a smooth elliptic subsolution.
For general cost-generated equations, boundary $C^{2,\alpha}$ estimates
under the strong Ma-Trudinger-Wang condition A3s were obtained in
\cite{HuangJiangLiu2015}, while interior H\"older and Sobolev estimates were developed in
\cite{LiuTrudingerWang2010,LiuTrudingerWang2015}.
For the second boundary value problem, global smoothness under the weak condition A3w was established in \cite{TrudingerWang2009}.
Just recently, Li, Wang, and Wang \cite{LiWangWang2026} proved global Sobolev and
H\"older estimates under A3w for the second boundary value problem on mutually uniformly $c$- and $c^*$-convex source and target domains. Hong and Wang \cite{HongWang2026} obtained boundary H\"older
regularity under A3s with weaker regularity assumptions on the domains. 

For the Dirichlet problem \eqref{eq:problem}, the image $T_u(\Omega)$ is determined by the
solution and has no prescribed boundary geometry. 
Consequently, the techniques developed for the second boundary value problem do not apply directly.
Our approach establishes the geometry of local sections near the boundary and develops perturbation arguments adapted to
the Dirichlet data.

\subsection{The a priori class}

We now specify the data on which the constants in our estimates may depend. 
Let $\Omega\subset\R^n$ be a bounded domain with $C^3$ boundary. Choose bounded open
ambient sets $U,V$, with $\overline\Omega\subset U$, and let $c$
be a fixed smooth cost on a neighbourhood of $\overline U\times\overline V$.
We impose the following assumptions; see \cite{MTW05,TrudingerWang2009}.
\begin{enumerate}
\item \textbf{(A1)} The twist maps $y\mapsto c_x(x,y)$ and $x\mapsto c_y(x,y)$
	are injective on their respective domains for each fixed $x$ and $y$,
	respectively. \\
\textbf{(A2)} The mixed Hessian $c_{xy}$ is
invertible, and both $c_{xy}$ and its inverse are uniformly bounded.
 
\item \textbf{(A3w)} The weak regularity condition holds:
\begin{equation*}
 A_{ij,p_kp_l}(x,p)\xi_i\xi_j\eta_k\eta_l\ge0
 \qquad \forall\, \xi\cdot\eta=0.
\end{equation*}

\item For every focus $y$, the image $c_y(\Omega,y)$ is
uniformly convex. The transformed boundaries have uniform $C^3$
bounds and a common modulus of continuity of their third derivatives.
The boundary datum $\varphi$, expressed in the corresponding coordinates, satisfies uniform $C^3$ bounds and a common modulus of continuity of its third derivatives.

\item A fixed compact set $V_0\Subset V$ contains
$T_u(\overline\Omega)$, where $T_u(x)=Y(x,Du(x))$. The region contains
$\overline\Omega\times V_0$ with a uniform neighbourhood. For every
$x\in\overline\Omega$ and $y_0,y_1\in V_0$, the segment
joining $c_x(x,y_0)$ to $c_x(x,y_1)$, together with a uniform neighbourhood
of that segment, lies in the inverse-twist domain. The cost derivatives,
inverse maps, and source convexity are uniformly controlled there.

More precisely, this control holds on the entire product
$\overline\Omega\times\{y\}$ for every focus $y$ arising from
such a segment and for every allowed small variation of that
focus. For all such foci, the images $c_y(\Omega,y)$ satisfy
the same uniform convexity and boundary regularity bounds.
These requirements also apply to the small variations of
supporting foci used below.

\item The solution $u\in C^4(\Omega)\cap C^2(\overline\Omega)$
is elliptic and satisfies
$\|u\|_{C^1(\overline\Omega)}\le M$. The density satisfies
$0<\lambda\le f\le\Lambda$.
\end{enumerate}

We denote this collection of a priori data by $\mathcal D$.
The collection $\mathcal D$ contains no a priori bound on $D^2u$ or  higher derivatives of $u$.

\subsection{Main results}

We first state the boundary localisation theorem, which provides
the section geometry needed for the global Sobolev estimate.
The global H\"older estimate combines this localisation with
a quantitative bound for the matrix term evaluated along the
solution. In all three statements, $\mathcal D$ denotes the
a priori data specified above.

\begin{theorem}[Boundary localisation]\label{thm:1}
Assume the a priori class in Section~1.1, let $n\ge2$, and suppose
only $\lambda\le f\le\Lambda$ for the density. For any
$z\in\partial\Omega$, subtract the tangent cost support and
use canonical source coordinates with $z=0$ and inward normal
$e_n$. Let $v$ be the resulting function, $D$ the transformed
source, $K_h=\{v<h\}$ and $d_h=\sup_{K_h}q_n$.
There are uniform $h_0,c,C>0$ depending only on $n,\mathcal D,\lambda,\Lambda$ and slides
\[
 T_h(q',q_n)=(q'-\nu_hq_n,q_n),\qquad
 |\nu_h|\le C(1+|\log h|),
\]
such that, for $0<h\le h_0$,
\begin{equation}\label{loc:1}
 B_{c\sqrt h}\cap T_hD\subset T_hK_h
 \subset B_{C\sqrt h}\cap T_hD,
 \qquad c\sqrt h\le d_h\le C\sqrt h.
\end{equation}

\end{theorem}

\begin{theorem}[Global Sobolev estimate]\label{thm:2}
Assume the a priori class in Section~1.1, let $n\ge2$, and suppose
that $f$ has a prescribed modulus of continuity on
$\overline\Omega$. For every $1<p<\infty$ there is a constant
$C_p$ depending only on $n,p,\mathcal D,\lambda,\Lambda$
and the modulus of continuity of $f$ such that every solution in this class satisfies
\[
 \|u\|_{W^{2,p}(\Omega)}\le C_p.
\]
 The constant is independent of
the Hessian and higher derivatives of the solution.
\end{theorem}

\begin{theorem}[Global H\"older estimate]\label{thm:3}
Assume the a priori class, let $n\ge2$, $0<\alpha<1$.
If $f\in C^\alpha(\overline\Omega)$, then
\[
 \|u\|_{C^{2,\alpha}(\overline\Omega)}\le C_\alpha,
\]
where $C_\alpha$ depends only on
$n,\alpha,\mathcal D,\lambda,\Lambda$ and
$\|f\|_{C^\alpha(\overline\Omega)}$.
\end{theorem}

For Theorem~\ref{thm:1}, we derive quadratic boundary separation and
two-sided section volume bounds from A3w and the Dirichlet trace.
These yield a boundary gradient modulus and compactness after
normalisation. Convex-envelope approximation and boundary comparison
identify the limiting problem, to which Savin's height-improvement
result \cite{Savin2013Boundary} applies. We then construct the
localising slides and obtain their logarithmic bound by comparison
across scales.

For Theorem~\ref{thm:2}, we combine the interior compactness and
section iteration of \cite{LiWangWang2026} with directional quadratic
bounds proved by smooth approximation. This gives the normalisation
needed for the perturbative Sobolev estimate of \cite{ChenFigalli2017}.
Near the boundary, Theorem~\ref{thm:1} controls maximal interior
sections and their John matrices. After normalisation and restriction
of the transport to a target ball and its inverse image, the same
perturbative estimate applies. Summing over dyadic boundary layers
gives the global $W^{2,p}$ bound.

For Theorem~\ref{thm:3}, the boundary gradient modulus makes the
canonical matrix along the solution a higher order perturbation.
Comparison with convex constant-density models and the boundary
estimate of \cite{Savin2013Boundary} yield a quadratic approximation,
which we iterate to obtain boundary expansions with exponent $\alpha$.
Theorem~\ref{thm:2} and Morrey's inequality then control the matrix
variations. Interior quadratic approximation and a covering argument
give the global H\"older estimate for the Hessian.

The paper is organized as follows. Section~2 establishes the preliminary
section estimates and convex-envelope approximation. Section~3 proves
Theorem~\ref{thm:1}. Section~4 develops the interior normalisation and
boundary-layer estimates and proves Theorem~\ref{thm:2}. Section~5
establishes the comparison and quadratic approximation estimates and
proves Theorem~\ref{thm:3}.

\vskip5pt
\noindent{\bf AI disclosure:}
The global $W^{2,p}$ regularity under A3s is a natural continuation of our earlier work \cite{HuangJiangLiu2015}.
We used ChatGPT (OpenAI) to assist in extending the global $W^{2,p}$ and $C^{2,\alpha}$
regularity results to the A3w setting and in preparing this
manuscript. The authors have reviewed all AI-assisted content
and take full responsibility for the mathematical results
and the final manuscript.

\section{Preliminary section estimates}

We write $W[u]=D^2u-A(x,Du)$. Unless otherwise indicated, constants
depend only on $n,\mathcal D,\lambda,\Lambda$. All sections below are
intersected with the source domain. Suprema over a section can equivalently
be taken over its closure.

\subsection{Coordinates and supports}
Fix a base pair $(x_0,y_0)$ and subtract
$m(x)=u(x_0)+c(x,y_0)-c(x_0,y_0)$. In the coordinates
$q=c_y(x,y_0)-c_y(x_0,y_0)$, write $x=X(q)$ and $v=(u-m)\circ X$.
If $J=DX$, direct differentiation gives
\begin{equation}\label{eq:1}
 D^2v-\widehat A(q,Dv)=J^TW[u]J,\qquad
 \widehat A(q,0)=\widehat A_p(q,0)=0,\qquad
 \widetilde f=(\det J)^2f\circ X.
\end{equation}
Indeed the terms containing $D^2X$ cancel in the first identity; differentiating
$c_y(X(q),y_0)=q+\mathrm{constant}$ gives the second canonical identity.
The bounds for $J,J^{-1}$ imply $0<\lambda_0\le\widetilde f\le\Lambda_0$.
Since the paired vectors transform as $J\xi,J^{-T}\eta$,
whose pairing is $\xi\cdot\eta$,  A3w is preserved. We omit the hat when working in these
coordinates. The phase assumptions ensure that all segments used in Taylor
expansions below are controlled.

\begin{lemma}\label{lem:1}
In canonical coordinates, for a unit vector $\xi$,
\begin{equation}\label{eq:2}
 |A(q,p)|\le C|p|^2,\qquad A(q,p)[\xi,\xi]\ge-C|p||\xi\cdot p|.
\end{equation}
Every tangent cost function is a global support with a unique contact.
Support-subtracted sections are convex in their canonical coordinates,
and $T_u$ is injective.
\end{lemma}
\begin{proof}
Taylor's formula gives
$A(q,p)=\int_0^1(1-t)A_{pp}(q,tp)[p,p]\,dt$.
Decompose $p=p_\perp+(p\cdot\xi)\xi$. A3w makes the term quadratic
in $p_\perp$ nonnegative after contraction against $\xi\otimes\xi$;
the remaining terms are bounded by $C|p||p\cdot\xi|$.
Consequently a straight restriction $F(t)=v(q_0+t\zeta)$ satisfies
$F''\ge-K|F'|$, with uniform $K$ on bounded chords. The same inequality
holds for a cost mountain relative to the base focus, since its augmented
Hessian is zero.

Scalar comparison shows that $F'$ cannot pass from nonnegative to negative.
On an interval where $F'\ge0$, $e^{Kt}F'$ is nondecreasing. In particular,
on $[0,1]$,
\begin{equation}\label{eq:3}
 F(t)-F(0)\le M_0t[F(1)-F(0)]^+,\qquad
 F(1)-F(0)\ge e^{-K}F'(0)\quad\hbox{if }F'(0)\ge0,
\end{equation}
where $M_0=e^K$ is sufficient. For the first inequality, if $F(t)>F(0)$,
let $s$ be the last preceding point with $F(s)=F(0)$. Integrating the
weighted monotonicity gives
$(F(t)-F(s))/(F(1)-F(s))\le e^K(t-s)/(1-s)\le e^Kt$.
The other case is immediate.

The canonical source is convex. Along a ray from the base contact,
$F'(0)=0$, so $v\ge0$. In the interior, ellipticity makes $F''>0$
whenever $F'=0$. Thus $F$ is strictly increasing away from the contact.
For a boundary base point every open chord enters the domain, and the same
argument applies. On any chord a rise followed by a fall is impossible;
this proves convexity of all sublevels. If two points have the same focus,
their global support inequalities force their additive constants to agree.
Uniqueness of contact then proves injectivity.
\end{proof}

\subsection{Volume and boundary separation}
Let $K_h=\{v<h\}$, where $v\ge0$ is a support remainder. Then there exists the John ellipsoid such that $a+LB_1\subset K_h\subset a+nLB_1$ . 

\begin{lemma}\label{lem:2}
For every such section, $|K_h|\le Ch^{n/2}$.
\end{lemma}
\begin{proof}
Put $E=a+(2M_0)^{-1}LB_1$. For $q\in E$, let $m_q$ be the
mountain supporting $v$ at $q$. It is at most $h$ on $a+LB_1$ and
$m_q(q)=v(q)\ge0$. Applying \eqref{eq:3} between $a$ and
$a+2M_0(q-a)$ gives $m_q(a)\ge-h$. Applying its second inequality
on $a\pm tLe$ in the direction with nonnegative initial derivative gives
$|L^TDm_q(a)|\le Ch$ for every unit $e$.
The fixed-source twist $y\mapsto c_q(a,y)$ has uniformly bounded inverse
Jacobian. Therefore $|T(E)|\le Ch^n/|\det L|$. On the other hand,
injectivity and the area formula give
\[
 |T(E)|=\int_E\frac{\widetilde f(q)}{|\det c_{qy}(q,T(q))|}\,dq
 \ge c|E|\ge c|\det L|.
\]
John's outer inclusion now proves the assertion. Only $E$ is used in
the area formula, so the argument also applies to boundary sections.
\end{proof}

Fix a boundary base point and rotate coordinates so that the inward normal
is $e_n$. Uniform source convexity and the boundary charts give
\begin{equation}\label{eq:4}
 D\subset\{q_n\ge\kappa|q|^2\},\qquad
 q_n=\rho(q'),\quad \rho(0)=D\rho(0)=0,\quad
 \kappa|q'|^2\le\rho(q')\le C|q'|^2.
\end{equation}
Write $d=d_h=\sup_{K_h}q_n$, $R=R_h=\sup_{K_h}|q|$ and
$g(q')=v(q',\rho(q'))$. The family of traces $g$ has uniform $C^3$
bounds and a common third-derivative modulus: on the boundary $u=\varphi$,
and the only solution-dependent parameters are the base point and focus.
In particular $0\le g(q')\le C|q'|^2$, and a fixed small $a>0$ gives
a boundary footprint $|q'|\le a\sqrt h$ with $g\le h/2$.

\begin{lemma}\label{lem:3}
For uniform $\mu,r_1,h_0>0$,
\begin{equation}\label{eq:5}
 \mu|q'|^2\le g(q')\le\mu^{-1}|q'|^2\quad(|q'|<r_1),
 \qquad ch\le d\le C\sqrt h,\quad a\sqrt h\le R,\quad R^2\le Cd.
\end{equation}
\end{lemma}
\begin{proof}
The midpoint of opposite footprint points gives $d\ge ch$, and
\eqref{eq:4} gives $R^2\le Cd$. If $P=(b,d)\in\overline K_h$,
choose signs $\sigma_i$ with $\sigma_i b_i\le0$. The simplex with
vertices $0,P,(\sigma_i a\sqrt h e_i,\rho(\sigma_i a\sqrt h e_i))$,
$1\le i<n$, has volume at least $ch^{(n-1)/2}d$ and lies in
$\overline K_h$. Lemma~\ref{lem:2} gives $d\le C\sqrt h$.
Thus $R\le Ch^{1/4}$ uniformly, and small sections lie in the chart.

If the smallest eigenvalue of $D^2g(0)$ had no positive uniform lower
bound, tangential rotation and compactness in $C^3$ would give
$g_j\to g_*$ with $D_{11}g_*(0)=0$. Nonnegativity implies
$D_{111}g_*(0)=0$, and hence $g_*(te_1)=o(|t|^3)$.
Fix a large $L$ and then a sufficiently small fixed $h$, and set
$r=Lh^{1/3}$. For large $j$, the points with tangential coordinates
$\pm re_1$ and $a\sqrt h e_i$, $2\le i<n$, have boundary value below
$h/2$. Together with zero they form a simplex of volume at least
\[
 \frac1{n!}r[\rho_j(re_1)+\rho_j(-re_1)](a\sqrt h)^{n-2}
 \ge cL^3h^{n/2}.
\]
This contradicts Lemma~\ref{lem:2}. For $n=2$ the product in the other
directions is empty. The uniform Taylor estimates now prove \eqref{eq:5}.
\end{proof}

\begin{lemma}\label{lem:4}
For small boundary sections,
\begin{equation}\label{eq:6}
 ch^{n/2}\le|K_h|\le Ch^{n/2},\qquad d_h\ge ch^{n/(n+1)}.
\end{equation}
If $q_h+L_hB_1$ is a John ellipsoid, then
\begin{equation}\label{eq:7}
 \|L_h\|\le R_h,\quad \|L_h^Te_n\|\le d_h/2,\quad
 \|L_h^{-1}\|\le\frac{CR_h}{d_h\sqrt h},\quad
 h\|L_h\|^2\|L_h^{-1}\|^2+h\|L_h^{-1}\|\le C.
\end{equation}
\end{lemma}
\begin{proof}
Put $L_0(z',z_n)=(\sqrt h z'+bz_n,dz_n)$, where $P=(b,d)$ is a
highest point. $L_0^{-1}P=e_n$,  and a footprint
point becomes
\[
 \left(z-\frac{b\rho(\sqrt h z)}{d\sqrt h},\frac{\rho(\sqrt h z)}d\right)
\] with the slides $L_0^{-1}$. 
By $|b|^2\le Cd$ and $d\ge ch$, this differs from $(z,0)$ by
$O(|z|^2)$. A simplex formed by $e_n$ and fixed sufficiently small
tangential points therefore contains a uniform ball. Its pullback lies
in $\overline K_h$, so John's inclusions give $\|L_h^{-1}L_0\|\le C$.
The width bounds and
$\|L_0^{-1}\|\le C(h^{-1/2}+|b|/(d\sqrt h)+d^{-1})$
prove the first three estimates in \eqref{eq:7}. The last follows from
$R^2\le Cd$, $d\ge ch$.

Since $0\in\overline K_h$, $|L_h^{-1}q|\le2n$ on $K_h$.
Consider
\[
 w(q)=\varepsilon h|L_h^{-1}q|^2+\varepsilon hq_n/d.
\]
On the level boundary $w\le\varepsilon(4n^2+1)h\le h$. On the
physical boundary of the section, \eqref{eq:5} gives $|q'|^2\le Ch$, and
\[
 h|L_h^{-1}q|^2\le C\left(|q'|^2+
 \frac{|b|^2+h}{d^2}\rho(q')^2\right)\le C|q'|^2,
 \qquad h\rho(q')/d\le C|q'|^2.
\]
Thus $w\le v$ on the whole boundary for fixed small $\varepsilon$.
In coordinates $q=L_h\xi$, the normalised matrix and density are
\begin{equation}\label{eq:8}
 A_h(\xi,p)=h^{-1}L_h^TA(L_h\xi,hL_h^{-T}p)L_h,
 \qquad f_h=\frac{(\det L_h)^2}{h^n}\widetilde f(L_h\xi).
\end{equation}
The normalised barrier has Hessian $2\varepsilon I$ and gradient
bounded by $C\varepsilon$. Its physical gradient is admissible by
\eqref{eq:7} and the phase margin. Moreover, by \eqref{eq:2} and
\eqref{eq:7}, we have $|A_h(\xi,D(w/h))|\le C\varepsilon^2$. Its augmented Hessian is at least $\varepsilon I$.
If $(\det L_h)^2/h^n<\varepsilon^n/\Lambda_0$, strict comparison
would give $w\le v$ in $K_h$. Along the inward normal this implies
$0=D_nv(0)\ge D_nw(0)=\varepsilon h/d>0$. This contradiction proves
the lower volume bound. Finally, \eqref{eq:4} places $K_h$ in a cap of
volume at most $Cd^{(n+1)/2}$, proving the lower bound for $d$.
\end{proof}

\subsection{Gradient bounds and convex approximation}
\begin{lemma}\label{lem:5}
Put $\alpha_0=(2n+1)^{-1}$. Then 
\begin{equation}\label{eq:9}
 0\le u(x)-m_z(x)\le C|x-z|^{1+\alpha_0},\qquad
 |Du(x)-Du(z)|\le C|x-z|^{\alpha_0},
\end{equation}
uniformly for $z\in\partial\Omega$ and
$x\in\overline\Omega$ near $z$. For boundary John matrices,
\begin{equation}\label{eq:10}
 \delta_h:=h\|L_h\|^2\|L_h^{-1}\|^2\le Ch^{1/(4n+1)},\qquad
 h\|L_h^{-1}\|\le Ch^{1/(8n+2)}.
\end{equation}
Thus every fixed normalised gradient ball is admissible for small $h$,
and  $|A_h|+|D_pA_h|+|D^2_{pp}A_h|\to0$ uniformly on the gradient ball.
\end{lemma}
\begin{proof}
Join $P=(b,d)$ to the footprint point with tangential coordinate
$-a\sqrt h b/|b|$. Their segment crosses the normal axis at height
at least
$a\sqrt h d/(|b|+a\sqrt h)\ge c\sqrt{hd}\ge ch^\gamma$,
where $\gamma=(2n+1)/(2n+2)$. If $b=0$, we consider $ P$ itself.
Choose $s=ch^\gamma$ so that $se_n\in\overline K_h$.
For $|q'|\le s/4$, the midpoint of $se_n$ and the footprint point
$(2q',\rho(2q'))$ lies above height $s/2$. The vertical segment joining
this midpoint to $(q',\rho(q'))$ lies in $\overline K_h$. Hence
$D\cap B_{s/4}\subset\overline K_h$. Inverting the relation between
$s$ and $h$ proves the first assertion of \eqref{eq:9}.

Ellipticity and the phase bounds give $D^2(u-m_z)\ge-CI$.
On a ball $B_s(x)\subset\Omega$ with $|x-z|\le C_1s$, two-sided
semiconvex difference quotients and the growth bound give
$|Du(x)-Dm_z(x)|\le Cs^{\alpha_0}$. First, we apply this with a nearest
boundary point to $x$. For boundary points $z,w$ with  $|z-w|=r$,
apply it to both supports at $z+4r\nu_z$, using a ball of radius $r$
in the uniform tubular neighbourhood. This gives the boundary-pair
gradient bound. Combining the two estimates proves \eqref{eq:9}.
In canonical coordinates it implies $|Dv(q)|\le C|q|^{\alpha_0}$.

The normal axis leaves $K_h$ at a point of height at most $d$, so
$h\le CR^{\alpha_0}d$. Integrating vertically from a point of $K_h$
to the boundary and using \eqref{eq:5} gives
$R^2\le C(h+R^{\alpha_0}d+d^2)$. Consequently
\[
 \begin{split}
 &h/d+R^2/d\le CR^{\alpha_0},\qquad
 \delta_h\le C(R^2/d)^2\le CR^{2\alpha_0},\\
 &h\|L_h^{-1}\|\le C\sqrt{(h/d)(R^2/d)}\le CR^{\alpha_0}.
 \end{split}
\]
Since $d\le C\sqrt h$, $R\le Ch^{1/(4-2\alpha_0)}$.
This proves \eqref{eq:10}. Taylor's formula in \eqref{eq:8} gives
$|A_h(\xi,p)|\le C\delta_h|p|^2$, $|D_pA_h|\le C\delta_h|p|$,
and $|D^2_{pp}A_h|\le C\delta_h$. Admissibility follows because
$hL_h^{-T}p\to0$ uniformly on each fixed ball.
\end{proof}

\begin{lemma}\label{lem:6}
Let $U(\xi)=h^{-1}v(L\xi)$ on a convex normalised section $O$ of
bounded diameter, with $0\le U\le1$. Put $\kappa=C\|L\|\sup|Dv|$.
Its lower convex envelope $\Gamma$ satisfies
\begin{equation}\label{eq:11}
 U_{ee}\ge-\kappa|U_e|,\qquad 0\le U-\Gamma\le C_n\kappa\DIAM O.
\end{equation}
On compact subsets at distance $r>0$ from the boundary,
$|DU|\le C_r$ and $D^2U\ge-C_r\kappa I$.
\end{lemma}
\begin{proof}
The directional inequality transforms as
\[
 U_{ee}=h^{-1}v_{Le,Le}\ge-C|Le||Dv||U_e|.
\]
A restriction to a segment is monotone or decreases then increases,
so its total variation is at most $2$. The Dirichlet Green function on
a segment of length $l$ is bounded by $l/4$. Integrating
$F''\ge-\kappa|F'|$ therefore bounds the error above its chord by
$\kappa l/2$. Every convex-envelope value uses a convex combination
of at most $n+1$ points. Successive applications of this chord bound
give \eqref{eq:11}. Integrating a positive directional derivative on
a segment of length $r/2$  with \eqref{eq:3}, it is bounded by $C_r$;
the opposite direction bounds a negative derivative. Substitution
in the directional inequality gives the Hessian bound.
\end{proof}

\section{Boundary localisation}

The volume estimate of Section~2 leaves open how far $K_h$ extends
in the normal direction. We apply Savin's height improvement to
limits of normalised sections to show that $d_h$ is comparable to
$\sqrt h$. We then construct the slides in Theorem~\ref{thm:1}
and compare them at successive dyadic heights to obtain the
logarithmic bound.

\subsection{Normalised sections and their limits}
After a rotation in the tangent plane we can choose a John matrix
\[
 L_h(\xi',\xi_n)=(D_h\xi'+\eta_h\xi_n,\tau_h\xi_n),
 \qquad D_h=\DIAG(\tau_{1,h},\ldots,\tau_{n-1,h}).
\]
This follows by a block Cholesky factorization of its shape matrix.
Set $O_h=L_h^{-1}K_h$, $U_h=h^{-1}v(L_h\xi)$,
$a_{i,h}=\tau_{i,h}/\sqrt h$, and $Q_h(\xi')=\sum_{i<n}a_{i,h}^2\xi_i^2$.
John's inclusions, widths, and \eqref{eq:6} imply
\begin{equation}\label{eq:12}
 B_1(\xi_h^*)\subset O_h\subset B_C\cap\{\xi_n>0\},\quad
 \tau_h\approx d_h,\quad \det L_h\approx h^{n/2},\quad a_{i,h}\ge c.
\end{equation}
Here and below $\approx$ denotes two inequalities with uniform constants.
To verify the last assertion and the boundary geometry, put
$\nu_h^J=\eta_h/\tau_h$. On the physical footprint the map
$P_h(q')=q'-\nu_h^J\rho(q')$ satisfies
\[
 \|DP_h-I\|\le CR_h\sqrt h/d_h
 =C[(R_h^2/d_h)(h/d_h)]^{1/2}\longrightarrow0.
\]
It is injective with a controlled inverse on a disk containing the footprint.
Using \eqref{eq:5}, its image contains $B_{c\sqrt h}$ and is contained
in $B_{C\sqrt h}$. Thus, for $G_h=L_h^{-1}(\partial D\cap\overline K_h)$,
\begin{equation}\label{eq:13}
 \begin{split}
 &G_h\subset\{0\le\xi_n\le Ch/d_h\},\qquad
 \{Q_h\le c\}\subset\pi(G_h)\subset\{Q_h\le C\},\\
 &U_h=1\ \hbox{on }\partial O_h\setminus G_h,\qquad
 cQ_h\le U_h\le\min(1,CQ_h)\ \hbox{on }G_h.
 \end{split}
\end{equation}
The set $G_h$ is a graph; its projection lies in $B_C$, proving
$\tau_{i,h}\ge c\sqrt h$. Moreover,
\begin{equation}\label{eq:14}
 ct^{n/2}\le|\{U_h<t\}|\le Ct^{n/2}\quad(0<t\le1),\qquad
 \lambda_*\le\det(D^2U_h-A_h)\le\Lambda_*.
\end{equation}
These follow from $\{U_h<t\}=L_h^{-1}K_{th}$ and \eqref{eq:8}.
Lemma~\ref{lem:6} gives convex $\Gamma_h$ with
$0\le\Gamma_h\le U_h\le1$ and $\|U_h-\Gamma_h\|_\infty\to0$.

We use the boundary convention:
the trace of a convex function is specified by its closed interior epigraph, see \cite[Definitions 2.1 and 2.5]{Savin2013Boundary}.
For a lower semicontinuous boundary datum $\phi$,  $\phi^*$ stands for
the trace of the closed convex hull of its epigraph.

\begin{lemma}\label{lem:7}
Every sequence $h_j$ decreasingly converging to 0 has a subsequence with a full-dimensional
convex limit domain $O$, a convex limit $V$, and a boundary datum
$\phi:\partial O\to[0,1]$, such that $U_j,\Gamma_j\to V$ locally
uniformly in $O$, the closed epigraphs of $U_j$ and $\Gamma_j$
converge to the closed epigraph of $V$, and
\begin{equation}\label{eq:15}
 \begin{split}
 &V(0)=0,\quad 0\le V\le1,\quad
 \lambda_*\le\det D^2V\le\Lambda_*,\quad V|_{\partial O}=\phi^*,\\
 &ct^{n/2}\le|\{V<t\}|\le Ct^{n/2}\quad(0<t\le1).
 \end{split}
\end{equation}
The raw boundary epigraphs of $U_j$ converge to that of $\phi$.
\end{lemma}
\begin{proof}
Truncate the epigraphs at height two. Blaschke compactness applies
to the closed domains and the closed convex epigraphs of $\Gamma_j$. The truncated
closed epigraphs of $U_j$ have the same limit, since their Hausdorff
distance is at most $\|U_j-\Gamma_j\|_\infty$. Extract a further
subsequence for the truncated boundary epigraphs. Each limiting boundary fiber is upward closed and
contains height one. This defines $\phi$. Convexity and boundedness
give local uniform convergence of $\Gamma_j$, hence of $U_j$.
Joining any point of the limiting epigraph to a fixed interior point at
height $3/2$ shows that it is the closure of its interior restriction.
In particular its trace contains $(0,0)$.

Test the equation for $U_j$, not the envelopes. At an upper contact
with a smooth test $\psi$, ellipticity gives
$D^2\psi-A_j\ge D^2U_j-A_j>0$ and determinant at least $\lambda_*$.
The test gradients are bounded, so Lemma~\ref{lem:5} permits passage
to the limit. At a lower contact with a nonnegative test Hessian, a
determinant exceeding $\Lambda_*$ would make that Hessian positive
definite. For large $j$ the augmented test Hessian remains admissible,
contradicting the equation. This proves the viscosity determinant
inequalities, equivalently the Alexandrov inequalities for convex $V$.

The limiting raw boundary epigraph lies in the full convex epigraph,
so $V|_{\partial O}\le\phi^*$. For the reverse inequality choose an
affine $l$ and $\eta>0$ with $l+2\eta\le\phi$ on $\partial O$.
Then $U_j\ge l+\eta$ on $\partial O_j$ for large $j$, by boundary
epigraph convergence. Choose smooth uniformly convex outer domains
$O^\varepsilon\supset\supset O$ converging to $O$ and solve
$\det D^2z^\varepsilon=2\Lambda_*$, $z^\varepsilon=0$ on
$\partial O^\varepsilon$; the classical Dirichlet theory supplies these
solutions \cite{CNS1984}. For fixed $\varepsilon$, their gradients and
positive definite Hessians are controlled on a common compact
neighbourhood of $O_j$. Thus $l+\eta/2+z^\varepsilon$ is a strict
augmented subsolution for large $j$, and comparison gives
$V\ge l+\eta/2+z^\varepsilon$. The Alexandrov estimate gives
$-z^\varepsilon(\xi)\le C\operatorname{dist}(\xi,\partial O^\varepsilon)^{1/n}$.
Taking the boundary trace and then $\varepsilon$ decreasing monotonically to 
0
 yields
$V\ge l+\eta/2$ there. Affine graphs of this kind separate every
point below the closed convex hull of the raw epigraph from that hull.
Indeed strict separation cannot have zero vertical coefficient, since
the point projects into $\overline O$; upward closure fixes its sign.
This proves the reverse trace inequality.

For $0<t<1$, interpolation with $(0,0)$ shows $\{V<t\}$ has interior.
Every point of $\{V=t\}$ lies on its boundary by joining it to a point
of smaller value. The level set has measure zero. Domain convergence
and local uniform convergence therefore give almost-everywhere convergence
of the sublevel indicators. Apply dominated convergence to \eqref{eq:14}, then
the same subsequence works for every $t$ increasing monotonically to 
1.
\end{proof}

\subsection{Normal height and slides}
Suppose $h_j$ decreases monotonically to 
0 and $d_{h_j}/\sqrt{h_j}\to0$. Ordering the
tangential axes and taking a subsequence, $a_{i,h_j}\to a_i<\infty$
for $i\le k$, and $a_{i,h_j}\to\infty$ for $k<i<n$. By
\eqref{eq:12}, $k\le n-2$. Set
$E_k=\{\xi_{k+1}=\cdots=\xi_n=0\}$ and $Q=\sum_{i\le k}a_i^2\xi_i^2$.
The limit $G$ of the physical boundaries satisfies
\[
 \{Q\le c\}\cap E_k\subset G\subset\{Q\le C\}\cap E_k.
\]
The raw datum is one off $G$ and is at least $cQ$ on $G$.
For the upper bound at a point with small $Q$, choose points on $G_j$
with the first $k$ coordinates fixed and the other tangential coordinates
exactly zero. Equation~\eqref{eq:13} gives $\phi\le CQ$ in the limit.
For larger $Q$ use $\phi\le1$. This choice of approximating points is
necessary because a vanishing coordinate times a diverging weight need
not vanish. Hence \eqref{eq:15} places $V$ with uniform $\mu$ in the class
$\mathcal D_0^\mu(a_1,\ldots,a_k,\infty,\ldots,\infty)$, see \cite[Section 5]{Savin2013Boundary}.
For $k=0$ this means $G=\{0\}$.

\begin{lemma}\label{lem:8}
There are uniform $h_*,b_*>0$ and $0<\theta<1$ such that, for
$b(h)=d_h/\sqrt h$,
\[
 h\le h_*,\quad b(h)\le b_*\quad\Longrightarrow\quad
 b(th)>2b(h)\ \hbox{for some }t\in[\theta,1].
\]
Consequently $c\sqrt h\le d_h\le C\sqrt h$ at all small heights.
\end{lemma}
\begin{proof}
Otherwise there are $h_j$ decreasing monotonically to 0, $b_j(h_j)\to0$ with
$b_j(th_j)\le2b_j(h_j)$ for $1/j\le t\le1$.
In the preceding normalisation,
\[
 \frac{\sup_{\{U_j<t\}}\xi_n}{\sqrt t}
 =\frac{\sqrt{h_j}}{\tau_{h_j}}b_j(th_j)\le C.
\]
Local uniform convergence gives the same bound for every fixed sublevel
of $V$. But for
each $M$ a member of the displayed infinite-weight class has
$\sup_{\{V<t\}}\xi_n\ge M\sqrt t$ at some $t\in[c_M,1]$, see \cite[Proposition 5.2]{Savin2013Boundary}.
This contradiction proves amplification.

For each fixed solution, strict increase along rays from zero implies
$\overline K_h=\overline D\cap\{v\le h\}$. Compactness of maximizing
points gives right continuity of $d_h$; radial contraction of such points
gives left continuity. Thus $b$ is positive and continuous. Equation
\eqref{eq:3} also gives that, for $h_1\le h_2$,
\begin{equation}\label{eq:16}
 M_0^{-1}\sqrt{h_1/h_2}\le b(h_1)/b(h_2)\le\sqrt{h_2/h_1}.
\end{equation}
Fix a uniform reference height $h_0\le h_*$, where \eqref{eq:6}
gives $b(h_0)\ge b_{\rm ref}>0$, and put $\beta=\min(b_*,b_{\rm ref})$.
If $b(h)<\beta/2$, let $h_1>h$ be its first crossing of $\beta/2$.
Amplification at $h_1$ gives $b(th_1)>\beta$, whence $th_1<h$ and
$h>\theta h_1$. Equation~\eqref{eq:16} yields
$b(h)\ge\sqrt\theta\beta/(2M_0)$. The upper bound is \eqref{eq:5}.
\end{proof}

\begin{proof}[Proof of Theorem~\ref{thm:1}]
Choose $P_h=(p_h',d_h)\in\overline K_h$ and set $\nu_h=p_h'/d_h$.
For small $h$, $P_h$ is off the physical boundary, whose section height
is $O(h)$, and $v(P_h)=h$. The slide $T_h$ sends it to $(0,d_h)$.
Moreover
\[
 |\nu_h|\sqrt h\le CR_h\to0,\qquad
 |\nu_h|^2\sqrt h\le CR_h^2/d_h\to0.
\]
In a cylinder of radius and height $c_1\sqrt h$, the slid boundary
satisfies $\xi_n=\rho(\xi'+\nu_h\xi_n)$. The derivative of
$\xi_n-\rho(\xi'+\nu_h\xi_n)$ in $\xi_n$ is at least $1/2$.
Thus it is a single graph $\widetilde\rho_h$, and absorption in
$\widetilde\rho_h\le2C|\xi'|^2+2C|\nu_h|^2\widetilde\rho_h^2$
gives $0\le\widetilde\rho_h\le C|\xi'|^2$. Its trace is at most
$C|\xi'|^2$.

For $X\in B_{c\sqrt h}\cap T_hD$, extend the line from $(0,d_h)$
through $X$ until it meets the physical boundary at $Y$. If $c$ is
small, $X_n<d_h/2$, and the continuation to height zero enlarges the
tangential coordinate by at most two. Thus $Y$ belongs to the graph
just described and has value strictly below $h$. Quasiconvexity bounds
the value on $[Y,(0,d_h)]$ by $h$. Equality at an interior point would
be a local maximum, contrary to the strict directional inequality.
This proves the inner inclusion in \eqref{loc:1}.

That inclusion contains a ball of radius $c\sqrt h$ centered at height
comparable to $\sqrt h$. The convex hull of this ball and any
$X\in T_hK_h$ has volume at least $ch^{(n-1)/2}$ times the distance
of $X$ from the ball. Since $\det T_h=1$ and $|K_h|\le Ch^{n/2}$,
we obtain the outer inclusion. If $h/2\le s\le h$, then
$P_s\in\overline K_h$ gives
$|(\nu_s-\nu_h)d_s|\le C\sqrt h$, so $|\nu_s-\nu_h|\le C$.
The slide at the fixed reference height is uniformly bounded.
Summing over dyadic heights yields $|\nu_h|\le C(1+|\log h|)$.
The normal-height assertion is Lemma~\ref{lem:8}.
\end{proof}

\section{The global Sobolev estimate}

We now assume that $f$ has a prescribed modulus of continuity.
In the interior, separation from tangent cost supports allows us to
normalise the sections and apply the perturbative $W^{2,p}$ estimate.
Near the boundary, we use Theorem~\ref{thm:1} to relate the height of
a maximal interior section to the distance of its center from
$\partial\Omega$. The local $L^p$ estimate on smaller sections then
gives a summable bound on dyadic boundary layers.

\subsection{Interior contacts}
\begin{lemma}\label{lem:9}
For $K\Subset\Omega$, the gradients have a common modulus of continuity
on $K$. For every $r>0$ there is $\sigma_K(r)>0$ such that
\begin{equation}\label{eq:17}
 u(x)-u(x_0)-c(x,T(x_0))+c(x_0,T(x_0))\ge\sigma_K(r)
 \quad(x_0\in K,\ |x-x_0|\ge r).
\end{equation}
These moduli use only the a priori data and $K$.
\end{lemma}
\begin{proof}
Consider a sequence $u_j$ converging uniformly to $u$. Uniform
semiconvexity implies, at interior points,
\begin{equation}\label{eq:18}
 \partial u(x)=c_x(x,\partial^cu(x)),\qquad \partial^cu(x)\subset V_0.
\end{equation}
Here $\partial u$ denotes the semiconvex subdifferential and
$\partial^cu$ the foci of global supports. Indeed, for $p\in\partial u(x)$,
a paraboloid of slope $p$ and sufficiently negative curvature touches
$u$ strictly from below on a small ball. Minimum points for $u_j$ minus
this paraboloid approach $x$, and their gradients approach $p$. Compactness
of $V_0$ and Lemma~\ref{lem:1} give a global limiting support with that
slope. The reverse inclusion follows from the supporting inequality.

The dual function $v(y)=\max_{\overline\Omega}(c(\cdot,y)-u)$ is locally
semiconvex. At its differentiability points, twist makes its maximizer
unique. Thus $E\mapsto|\partial^cu(E)|$ is a measure in the interior.
The area formula for $T_j$ gives two-sided bounds by $|E|$. For an open
$G\Subset\Omega$, almost every $y\in\partial^cu(G)$ has its unique
maximizer in $G$, and all maximizers for $u_j$ eventually lie in $G$.
Hence $y\in T_j(G)$. For compact $F\Subset\Omega$, closedness of the
support graph gives $T_j(F)\subset N_\varepsilon(\partial^cu(F))$ for
large $j$. Fatou's lemma and outer regularity, followed by regularity
of the measures, yield
\begin{equation}\label{eq:19}
 \lambda_1|E|\le|\partial^cu(E)|\le\Lambda_1|E|
 \quad(E\Subset\Omega\hbox{ Borel}).
\end{equation}

By \eqref{eq:9}, the boundary gradients converge along a subsequence,
and the growth estimate identifies their limit with the boundary jet of
$u$. A limiting boundary tangent support has unique contact: its
approximating sections have uniform diameter at most $Ch^{1/4}$.
If an interior support also touched the boundary at $z$, its contact
set would be convex in its own canonical coordinates. This follows by
applying the directional inequality to $u_j$ minus that fixed mountain
and passing to the limit. The chord to its interior contact is transverse
to the boundary. The remainder has zero derivative on this chord and
on the tangent plane, so its full boundary gradient is zero. Twist then
identifies it with the boundary tangent support, a contradiction.

We give the local domain verification for the Alexandrov estimates.
At $x_*$, $P=\partial u(x_*)$ is compact and convex. By \eqref{eq:18}
and the phase margin it has a closed neighbourhood in the inverse-twist
chart. Choose a smooth uniformly convex $P'$ with
$P\Subset P'$ in this neighbourhood, and put $V'=Y(x_*,P')$.
Such $P'$ can be obtained by smoothing the squared distance to $P$,
adding a small positive quadratic, and choosing an intermediate regular
sublevel. For $x$ sufficiently near $x_*$,
$p\mapsto c_x(x,Y(x_*,p))$ is $C^2$-close to the identity, so its image
of $P'$ is uniformly convex. For small $\rho$, the sets
$c_y(B_\rho(x_*),y)$ are also uniformly convex for
$y\in\overline V'$, by rescaling to nonsingular linear maps.
Closedness of the support graph puts all foci near $x_*$ in $V'$.
The restricted normal map on this product equals the global normal map
at interior contacts: a local support gives a subgradient, and
\eqref{eq:18} and twist identify its global focus.

We may therefore use the local Alexandrov estimates
\cite[Theorems 6.2 and 6.11]{FigalliKimMcCann2013} on sufficiently
small compact sections in this product with \eqref{eq:19}. In the
present sign convention, there are two consequences below:
\begin{equation}\label{eq:20}
 |Q|^2\le C|\min_Q w|^n,\qquad
 |w(q)|^n\le C(1-t_q)^{1/2^{n-1}}|Q|^2,
\end{equation}
where $w$ is a support remainder zero on $\partial Q$, and $t_q$ is
its radial fraction from the John center to the boundary. The second
inequality is used for $t_q>1/(2n)$; points with smaller radial fraction
already satisfy the centering conclusion needed below. The coordinate
Jacobian factors are bounded on the local product.

If a contact set had more than one point, its canonical image would
be a compact convex set in the interior. Choose an exposed point
$q_*$, exposed by $e\cdot q=b$, and $d>0$ so small that its cap
$e\cdot q\ge b-2d$ lies in a ball whose fivefold dilate is contained
in the preceding local chart. Decrease $d$ to meet the diameter
thresholds in the two Alexandrov estimates. The cap contains a segment
whose $e$-coordinate ranges from $b-d$ to $b$. Tilt the focus
by $\varepsilon e$ and adjust the constant. The new mountain obeys
\[
 m_\varepsilon-m=\varepsilon(e\cdot q-b+d)+O(\varepsilon^2).
\]
The section $Q_\varepsilon=\{u\le m_\varepsilon\}$ lies in the cap
neighbourhood for small $\varepsilon$: any limit point is an old contact
with $e\cdot q\ge b-d$. The old contact with $e\cdot q=b-d/2$
and the point $q_*$ both belong to $Q_\varepsilon$ for small
$\varepsilon$, by the expansion above. In the canonical coordinates
of the tilted focus the section is convex, and the coordinate change
tends uniformly to the identity on the cap. Thus its width in the
$e$ direction is at least $d/3$, while the image of $q_*$ has distance
$o(1)$ from its forward supporting plane. By John's inclusions the
John center has distance at least $d/[3(n+1)]$ from that plane.
Radial interpolation from this center gives $1-t_{q_*}=o(1)$. Moreover
\[
 (m_\varepsilon-u)(q_*)=\varepsilon d+O(\varepsilon^2),\qquad
 \max_{Q_\varepsilon}(m_\varepsilon-u)=\varepsilon d+o(\varepsilon).
\]
The latter follows by taking limit points of maximizers in the old contact
set. The two inequalities in \eqref{eq:20} contradict their ratio tending
to one. This proves strict contact.

If $\partial u(x)$ were nontrivial, choose an exposed slope and subtract
its cost support in canonical coordinates. After a linear change there
are a unit $e$ and $b>0$ such that the new remainder satisfies
$w(-te)=o(t)$ and $w(te)\ge bt-O(t^2)$. These follow from the support
function formula for directional derivatives of a semiconvex function.
Strict contact makes $Q_h=\{w\le h\}$ shrink to zero. Apply
\eqref{eq:20} at its minimum to $w-h$, then we get its John radial fraction is
at most $1-c$, so its two radii along $e$ and $-e$ are comparable.
The displayed one-sided growth contradicts this, thus every limit is
differentiable. Failure of \eqref{eq:17} would give a limiting support
with two contacts. Failure of gradient equicontinuity would give two
different limiting slopes at one interior point by semiconvexity,and 
both are impossible.
\end{proof}

\subsection{Interior normalisation}
We record the directional estimate needed for the interior compactness
argument. Its approximation is included to keep track of the domains
and constant densities in that argument.

\begin{lemma}\label{lem:10}
Let
$c_0(X,Y)=X\cdot Y+R(s,z)$, $X=(s,t)$, $Y=(\eta,z)$ in $\R^k\times\R^{n-k}$.  Suppose $R$ is finite, nonnegative,
separately convex, zero on both axes, and
$D_{s,\xi}^2D_{z,\zeta}^2R\ge0$ in distributions. Let finite convex
potentials $U,V$ satisfy $U+V\ge c_0$, and the equality holds on a plan of
constant positive marginals. Assume $U(0)=V(0)=0$, there exists uniform positive
gaps on compact annuli and bounded opposite coordinates on the plan
over each bounded ball in both directions. Then for uniform $b,C>0$,
\begin{equation}\label{eq:21}
 U(0,t)\le C|t|^2,\qquad V(\eta,0)\le C|\eta|^2
 \quad(|t|,|\eta|<b).
\end{equation}
The constants depend only  on these bounds and local upper bounds for $R,U,V$.
\end{lemma}
\begin{proof}
Restrict the plan to $B_a\times\R^n$, with $a$ large enough that its
target marginal remains $g_0\,dY$ on $B_4$. Its source marginal is
$f_0\mathbf1_{B_a}\,dX$ and both marginals have compact support.
Convolve $R$ to obtain smooth $R_j$. Its separate convexity and mixed
fourth-derivative sign persist. The cost $c_j=X\cdot Y+R_j(s,z)$ has
twists
\[
 (c_j)_X=(\eta+(R_j)_s,z),\ (c_j)_Y=(s,t+(R_j)_z) \ \text{and}\ \det(c_j)_{XY}=1.
\]
Its matrix is $A_j(s,p)=\DIAG((R_j)_{ss}(s,p_t),0)$, which is
matrix-convex in $p$.

For completeness, we choose mutually uniformly cost-convex smooth domains
\[
 D_j=\{|s|^2/S^2+|t|^2/T_j^2<1\}\ \text{and}\ 
 D_j^*=\{|\eta|^2/E_j^2+|z|^2/Z^2<1\}.
\]
Fix $S,Z$ large enough for the marginal supports, and then choose
$T_j,E_j$ large. For the source defining function, the corrected Hessian is
\[
 2\DIAG\left(S^{-2}I-T_j^{-2}\sum_\alpha t_\alpha(R_j)_{ssz_\alpha},
 T_j^{-2}I\right)>0
\]
if $T_j>2S^2\sup|(R_j)_{ssz}|$. To verify global convexity of its
shear image, we set $\tau=t+(R_j)_z$. The Hessian of the transformed
defining function on $(v,w)$ equals
\[
 \frac{2|v|^2}{S^2}+\frac{2|w-D_s(R_j)_zv|^2}{T_j^2}
 -\frac2{T_j^2}(\tau-(R_j)_z)\cdot D_s^2(R_j)_z[v,v].
\]
It is positive on a convex cylinder containing the image if we increase
$T_j$ so that
\[
 T_j>4S^2\sup|(R_j)_{ssz}|+2\sup|(R_j)_z|.
\]
Interchanging source
and target gives the same verification for $E_j$. These corrected
defining functions also supply the global convexity barriers needed
in the smooth second boundary value theory.

We convolve the given finite marginals and the convolved marginals coincide with $f_0,g_0$ on $B_3$.
We then add smooth background densities of equal mass $\delta_j$, supported 
outside $B_2$ and strictly positive on  $\overline B_2$, so that the total densities
are positive up to the boundaries. The parameter $\delta_j$ is chosen so that
$\delta_j(1+\sup_{D_j\times D_j^*}|c_j|)\to0$.
By the smooth theory \cite[Theorems 1.2 and 1.3]{TrudingerWang2009}, there exists optimal potentials $U_j,V_j$ and inverse diffeomorphisms. Moreover, these potentials are classically convex, since the cost $c_j$ is separately convex.

Every weak limit of their finite graph plans is optimal for the original
finite marginals. Indeed， the convolved masses have fixed compact
supports, the cost contribution of the backgrounds tends to zero, and
$c_j\to c_0$ uniformly on products of  compact sets. Independently convolving
the original plan and coupling the background densities gives a competitor with
convergent cost, proving optimality of the limit.
We normalise each potential by its value at zero, and they are locally
equi-Lipschitz on $B_2$. If the gradient of a convex potential is large at some point, convexity forces the gradient to remain at least half its size on a cone of fixed volume in a slightly
larger ball. On products of  compact sets, $(c_j)_X$ is bounded by separate
convexity, so this cone would be mapped outside a fixed ball containing all
convolved target supports. The resulting positive mass in the image contradicts the fact that the exterior
mass $\delta_j\to0$. The dual argument is identical.
Consider a subsequence along which the potentials converge locally to $\widehat U$ and  let $\pi$ be a
weak limit of the graph plans. At each $(X,Y)\in\SUPP\pi$ with
$X\in B_2$, pass to the limit in the supporting inequality along 
approximating  pairs of points in the graphs to obtain the supporting inequality for
$\widehat U$ on $B_2$. The restricted original plan has the same
marginals and the same optimal cost as $\pi$, and is supported on 
$\{U+V=c_0\}$, hence $\int_{X\times Y}(U+V-c_0)\,d\pi=0$. The integrand is
continuous and nonnegative, so equality holds everywhere on $\SUPP\pi$.
Thus $U$ satisfies the same supporting inequality as $\widehat U$, and the nonempty convex
subdifferential of $c_0(\cdot,Y)$ at $X$ is contained in both
$\partial\widehat U(X)$ and $\partial U(X)$.
At almost every  point in the source space, both functions are differentiable, so their gradients
agree. The normalisation condition then identifies them on $B_2$; similarly the dual limit
coincides with $V$. Thus the whole sequence converges locally uniformly.

Separate convexity together with the vanishing of R along the coordinate axes implies $0\le R(s,z)\le C|s||z|$
on products of sufficiently small fixed sets. For a  pair in the support of an optimal plan with $X$ small and $Y$ bounded, we have 
$0\le V(Y)\le c_0(X,Y)\le C|X||Y|$. The annular gap structure forces $Y\to0$,
and hence $U(X)=o(|X|)$; dually $V(Y)=o(|Y|)$. The convexity of the potentials then implies, via difference
quotients , that $DU_j,DV_j$ uniformly small on any fixed sufficiently
small ball. Furthermore, on this ball we have  $|(R_j)_s(s,z)|\le C(|z|+\varepsilon_j)$ and
$|(R_j)_z(s,z)|\le C(|s|+\varepsilon_j)$ there. The shear formulas
therefore map both normal images into $B_2$, and consequently both equations 
have constant right-hand sides on this fixed ball.

We subtract the one-variable gauges $R_j(s,0)-R_j(0,0)$ and
$R_j(0,z)-R_j(0,0)$ from the normalised potentials to obtain $u_j$ and $v_j$.
The corrected remainder is given by
$\widetilde R_j=R_j(s,z)-R_j(s,0)-R_j(0,z)+R_j(0,0)$.
The source matrix remains matrix-convex, is independent of $t$, and satisfies
$A(s,0)=0$. These gauge functions converge  locally uniformly to zero while preserving the 
derivatives in the linear directions. Futhermore, the annular gap structure yields a component
$G_j=\{u_j<\tau\}$ that is compactly contained within the constant-density ball,
containing a fixed ball $B_b\subset\{u_j<\tau/2\}$.

Fix a unit vector $a$ in the $t$ direction and omit the index $j$. For
$\mathcal Lw=W^{ij}(w_{ij}-A_{ij,p_l}w_l)$, differentiating  the
constant equation yields
\[
 \mathcal Lu_a=0,\ 
 \mathcal Lu_{aa}\ge\TR(W^{-1}W_aW^{-1}W_a),\ 
 \mathcal Lu\le n\ \text{and}\ \mathcal L(u_a^2)=2u_{aa}.
\]
Here,  the second inequality follows  matrix convexity, while the third follows from
$A-\sum_l A_{,p_l}p_l\le0$. The last identity uses $Du_a=Wa$.
Define  $M=\sup_{G_j}|u_a|>0$. At a maximum point of the auxiliary function 
$\log u_{aa}+\log(\tau-u)+u_a^2/(2M^2)$, we choose a constant
$W$-orthogonal frame whose first basis vector is $a$, and denote  $W_{ii}=\lambda_i$.
Then $\lambda_1=u_{aa}$ and $(W_a)_{i1}=u_{aai}$ since $Aa=0$.
After subtracting the squared logarithmic gradient, the trace term
reduces to $\sum_{i>1}u_{aai}^2/(\lambda_i u_{aa}^2)$.
For $i>1$, we have $u_{ai}=0$ and the first derivative test gives
$u_{aai}/u_{aa}=u_i/(\tau-u)$. These terms cancel the corresponding
contributions from $\mathcal L\log(\tau-u)$ yielding 
\[
 0\ge\frac{u_{aa}}{M^2}-\frac n{\tau-u}
       -\frac{u_a^2}{u_{aa}(\tau-u)^2}.
\]
Thus, $(\tau-u)u_{aa}\le C_nM^2$ holds at this maximum point. Since the exponential
factor lies  between 1 and $e^{1/2}$,  the same bound, up to a modified
constant, holds throughout $G_j$. The maximum must be interior because 
$u=\tau$ on the  boundary of $G_j$. Uniform bounds on the image sets provide a uniform control on $M$.
Integrating along the linear directions on $B_b$ and passing to the
limit yields  \eqref{eq:21}, as the linear terms tend to zero. Finally, interchanging
the variables establishes the corresponding dual assertion.
\end{proof}

\begin{proposition}\label{prop:1}
Let $K\Subset\Omega$ and $1<p<\infty$.
 If $f$ has a prescribed modulus of
continuity, then $\|D^2u\|_{L^p(K)}\le C_{K,p}$. In double canonical coordinates at contacts based at points of $K$,
let $L_h$ denote the John matrix of the source section, and set $P_h=hL_h^{-T}$.
Then for every $\varepsilon>0$, uniformly in $K$, there holds 
\begin{equation}\label{eq:22}
 \|L_h\|+\|P_h\|\le C_\varepsilon h^{1/2-\varepsilon},\qquad
 \|L_h^{-1}\|+\|P_h^{-1}\|\le C_\varepsilon h^{-1/2-\varepsilon}
\end{equation} 
for all sufficiently small $h$.
Moreover, every normalised potential limit is a positive definite quadratic form whose Hessian admits uniform upper and lower bounds.
\end{proposition}
\begin{proof}
Choose $K\Subset K_1\Subset\Omega$ and a uniform ball $B_{4R}(x_0)\Subset K_1$.
By Lemma~\ref{lem:9}, there exists $r>0$ so that, for every
$y\in B_{4r}(y_0)$ with $y_0=T(x_0)$, all maximizers of $c(\cdot,y)-u$
lie in $B_R(x_0)$. These maximizers are unique by Lemma~\ref{lem:1}.
Hence this target ball is filled by the actual map. Setting
$V=B_{2r}(y_0)$, $S=T^{-1}$ and $D=T^{-1}(V)$, differentiation
yields $c_{xy}DT=W[u]$, so that $T:D\to V$ is a diffeomorphism and
\begin{equation}\label{eq:23}
 \rho(x)=\frac{f(x)}{|\det c_{xy}(x,T(x))|},\qquad
 T_\#(\rho\mathbf1_D\,dx)=\mathbf1_V\,dy.
\end{equation}
The source contains a common ball about $x_0$, and
$\omega_\rho(t)\le C(\omega_f(t)+t+\omega_{Du}(t))$.

The dual remainder at $y_0$ has a strict-contact modulus. Indeed, consider the contact of
$v(y)=\max(c(\cdot,y)-u)$ at $x_0+ae$, where $e$ is the unit vector in the direction of 
$c_x(x_0,y)-c_x(x_0,y_0)$. The corresponding remainder is bounded below by
$a(c|y-y_0|-\omega_{Du}(a)-Ca)$; it therefore suffices to choose $a$ small, depending only on the
prescribed separation. Moreover, adding the two support inequalities
at $x_i=S(y_i), i=1,2$, yields
$\sigma_{K_1}(|x_1-x_2|)\le C|y_1-y_2|$.
Consequently, $S$ and $Dv$ also share a common local modulus of continuity.

The ambient normal-measure hypotheses on both sides require careful treatment.
On $U=B_{2R}(x_0)$, we define the extension
$\widetilde u(x)=\max_{\overline V}(c(x,y)-v(y))$.
The function  agrees with $u$ on $D$, and its transform on $V$ coincides with $v$.
If  $y\in V$ is an interior point that supports $\widetilde u$ at $x$, then differentiating the dual
inequality gives $Dv(y)=c_y(x,y)$, whence $x=S(y)$. Since boundary foci form a set of measure zero, the upper normal-measure bound holds on $U$, while the lower bound holds on $D$. On the other side, choose a common ball
$B= B_b(x_0)\Subset D$ with $T(B)\Subset V$, and define 
$\widetilde v(y)=\max_{\overline B}(c(x,y)-u(x))$, $y\in V$. By the inverse modulus, $\widetilde v$ equals $v$ on a fixed smaller ball about $y_0$. If $x\in B$ is an interior normal point of $\widetilde v$ at $y$, then necessarily $y=T(x)$;
all remaining normal values lie on $\partial B$. Hence
$|\partial^{c^*}\widetilde v(E)|\le C|E|$, with the reverse inequality holding on that smaller ball.  After reducing their radii, the ambient balls are mutually uniformly
cost-convex. By the virtue ofthe two contact moduli, small sections and their normal
images remain within the regions of agreement.

We now pass to double canonical coordinates, in which 
$c(q,0)=c(0,\eta)=0$, $c_q(0,\eta)=\eta$, $c_\eta(q,0)=q$.
Let $x=\chi(q)$ and $y=\psi(\eta)$ denote the inverse coordinate maps, then
the densities in \eqref{eq:23} become
\[
 \rho_c(q)=\rho(\chi(q))|\det D\chi(q)|,\qquad
 \rho_c^*(\eta)=|\det D\psi(\eta)|.
\]
Both densities admit uniform positive upper and lower bounds and share a common modulus
of continuity on the respective local domains. Denote $u,v$ by the corresponding
support-subtracted potentials, and let $S_h,S_h^*$ denote the sections of $u$ and $v$ at zero.
The verifications carried out above permit us to apply the local section estimates of  
\cite[Theorems 6.2, 6.11, 9.3 and Lemma 9.4]{FigalliKimMcCann2013}
and \cite[Proposition 2.1(ii)]{DePhilippisFigalli2013} on both sides.
The polar comparison and centering argument
\cite[Lemma 5.1 and Proposition 5.2]{LiWangWang2026} then yield
\[
 B_{1/C}\subset L_h^{-1}S_h\subset B_C,\qquad
 B_{1/C}\subset P_h^{-1}S_h^*\subset B_C,
 \qquad \|L_h\|+\|P_h\|\to0.
\]
In particular, the target John matrix can be written as $P_h C_h$, where
$\|C_h\|+\|C_h^{-1}\|\le C$. The volume bounds for the source section give
$c\le|\det L_h|/h^{n/2}\le C$. For $H\ge1$ with $Hh$ sufficiently
small, height comparison gives constants $0<\beta_1\le\beta_2$, such that
$C^{-1}H^{\beta_1}S_h\subset S_{Hh}\subset CH^{\beta_2}S_h$,
and the analogous inclusions hold for the dual sections.
The engulfing property then yields $C^{-1}u(q)\le v(\eta)\le Cu(q)$ on the graph of the map.
Consequently, both normalised domains exhaust $\R^n$,  the potentials remain
bounded on fixed balls, and the coordinates of the opposite graph are bounded there.
Introducing the variables $q=L_hX$ and $\eta=P_hY$, we divide both transported
measures by $h^{n/2}$. The resulting marginal densities are
\[
 F_h(X)=\frac{|\det L_h|}{h^{n/2}}\rho_c(L_hX),\qquad
 G_h(Y)=\frac{|\det P_h|}{h^{n/2}}\rho_c^*(P_hY).
\]
The two prefactors are reciprocal to each other, and their common two-sided bounds
remain valid as $h$ decreasing monitonically to 0. Along a subsequence, $F_h(0)$ and $G_h(0)$
converge to positive constants. The common continuity moduli, together with the uniform
shrinking of $L_h,P_h$, imply convergence to those constants uniformly
on each fixed ball.
These conditions constitute the local hypotheses of
\cite[Proposition 5.3]{LiWangWang2026}; in particular, they do not require any boundary
geometry of the actual image $T(\Omega)$.

This compactness result yields convex limits with annular gaps. The cost
limit is $X\cdot Y+R(X_{I_+},Y_{I_-})$.
The sets $I_+,I_-$ are the semi-axes whose
source semi-axes divided by $\sqrt h$ tend to infinity and
zero, respectively. The remainder possesses exactly the properties stated in Lemma~\ref{lem:10}.
Let $\Lambda_h=h^{-1}\|L_h\|^2\|P_h\|^2$.
If $\Lambda_{h_j}\ge\epsilon>0$, then both extreme sets are nonempty.
The constants in \eqref{eq:21} are uniform along these limits, as a consequence of the
annular gaps, local upper bounds and graph bounds established above. For
$\sqrt\theta<b$ in Lemma~\ref{lem:10}, testing the dual inequality with the vector
$Y_{I_+}$ of length $\sqrt\theta$ in the direction of $X_{I_+}$
gives $|X_{I_+}|\le C\sqrt\theta$ on $\{U\le\theta\}$.
The analogous bound holds for $|Y_{I_-}|$. The other semi-axes contribute terms of order
$o(\|L_{h_j}\|)$ and $o(\|P_{h_j}\|)$ to the physical lengths.
Returning to the inner ellipsoid inclusions at height $\theta h_j$, we obtain
$\Lambda_{\theta h_j}\le C^4\theta(1+o(1))\Lambda_{h_j}$.
Fix $\theta$ with $\sqrt\theta<b$ and $C^4\theta<1/4$, where $\theta$ is independent
of $\epsilon$. For each $\epsilon>0$, the preceding subsequence argument
gives a common $h_\epsilon>0$ such that $\Lambda_{\theta h}\le
\Lambda_h/2$ whenever $h\le h_\epsilon$ and $\Lambda_h\ge\epsilon$.
The nesting of both families of sections together with their centered ellipsoid
inclusions, yields $\Lambda_{th}\le C_\theta\Lambda_h$ for
$\theta\le t\le1$. Setting $h_k=\theta^kh_\epsilon$, we obtain
\[
 \Lambda_{h_{k+1}}\le\max\{\Lambda_{h_k}/2,C_\theta\epsilon\}.
\]
At the fixed height $h_\epsilon$, $\Lambda_{h_\epsilon}$ is uniformly
bounded by the diameter bounds for the two sections. Iterating this estimate and
then applying the nesting bound on $[h_{k+1},h_k]$, we obtain 
$\limsup_{h\to0}\Lambda_h\le C_\theta^2\epsilon$ uniformly.
Letting $\epsilon$ decreasing monotonically to 0, we conclude that$\Lambda_h\to0$.

The remaining bilinear-limit and section iteration argument of
\cite[Proposition 5.7]{LiWangWang2026} now applies: all the limiting convex potentials have constant Monge--Amp\`ere density and are quadratic,
and comparison of their sections at a fixed smaller height gives
\eqref{eq:22}. In particular, the normalised costs converge to $X\cdot Y$
in $C^2$ on fixed products. Upon choosing $\varepsilon<1/6$, their third derivatives are bounded by
$Ch^{-1}(\|L_h\|+\|P_h\|)^3=o(1)$.
We restrict to a convex dual section and its actual inverse image.
Centering, height comparison, and the engulfing property ensure that both are contained between fixed
inner and outer balls. Set $a=(G_h(0)/F_h(0))^{1/(2n)}$,
$X=a\widehat X$ and  $Y=a^{-1}\widehat Y$, which preserves the bilinear form $X\cdot Y$.
The two centering factors now coincide and  equal
$m=a^nF_h(0)=a^{-n}G_h(0)$. Dividing both measures by $m$, we obtain the
densities $F_h(a\widehat X)/F_h(0)$ and $G_h(a^{-1}\widehat Y)/G_h(0)$
on the restricted domains. Both densities converge uniformly to the constant 1. The bounds
on $a$ and $a^{-1}$ ensure that the fixed inner and outer balls are preserved up to a common
change of radii,and that the  convergence of the costs on fixed
products is maintained. Applying
\cite[Proposition 2.4 and Section 3.4]{ChenFigalli2017}, with  the minimizing
cost taken to be the negative of our cost, yields a uniform local $L^p$
Hessian bound, provided that $h$ is sufficiently small for the prescribed $p$.
At this fixed height, the section volume and diameter bounds yield uniform control of
$L_h$ and $L_h^{-1}$. The pullback along the smooth canonical coordinates introduces  only terms controlled by $\|u\|_{C^1}$. A finite covering argument for
$K$ then completes the proof of the asserted estimate.
\end{proof}

\subsection{Maximal interior sections}
Denote $\ell(s)=1+|\log s|$ for small $s>0$. For $x\in\Omega$, let
$m_x$ be its tangent support and set
\[
 w_x=u-m_x,\qquad H(x)=\min_{\partial\Omega}w_x,\qquad
 S_t(x)=\{w_x<t\}.
\]
The uniqueness of the contact point guarantees that $H>0$ and that $S_t(x)\Subset\Omega$ for $t<H$; moreover, 
the initial gradient bound gives $H\le C\operatorname{dist}(x,\partial\Omega)$.

\begin{proposition}\label{prop:2}
For sufficiently small $H=H(x)$, in double canonical source coordinates $Q$ at
$(x,T(x))$, any John matrix $L$ for $Q(S_{3H/4}(x))$ satisfies
\begin{equation}\label{eq:24}
 \begin{split}
 &H\approx\operatorname{dist}(x,\partial\Omega)^2,\qquad
 S_H(x)\subset\{\operatorname{dist}(\cdot,\partial\Omega)<C\sqrt H\},\\
 &|\det L|\approx H^{n/2},\qquad
 \|L\|\le C\sqrt H\ell(H),\qquad
 \|L^{-1}\|\le C\ell(H)/\sqrt H.
 \end{split}
\end{equation}
\end{proposition}
\begin{proof}
Choose a boundary contact point $z$ of $m_x+H$ and introduce double canonical
coordinates $(q,\eta)$ at $(z,T(z))$. The lift of the support function has parameter
$ae_n$, since the tangential derivatives agree at $z$. Set
$g_a(q)=c(q,ae_n)$ and $v=u-m_z$. Then
\begin{equation}\label{eq:25}
 S_H(x)=\{v<g_a\},\qquad H=\max_D(g_a-v)=(g_a-v)(q_x).
\end{equation}
We first show that $a>0$. In the source coordinates $Q_a=c_\eta(q,ae_n)$,
the domain has the same tangent plane at the origin, since $DQ_a(0)=I$.
Along the chord joining the origin to $Q_a(q_x)$, the function $w_x$ takes the  value $H$ and 0, with
initial derivative $-a(Q_a(q_x))_n$. If $a\le0$, the directional
inequality would force this derivative to remain nonnegative, a contradiction.

 Decreasing $c$ outside a fixed
small section, the normal-height upper bound of Theorem~\ref{thm:1} gives
$v(q)\ge cq_n^2$ throughout $D$. The inner inclusion then yields $v(te_n)\le Ct^2\ell(t)^2$.
Indeed, choose $s=C_1t^2\ell(t)^2$ so that $|T_s(te_n)|<c\sqrt s$.
For any fixed small admissible $a_0>0$, the source Taylor expansion gives
$g_a(te_n)=at+O(t^2)$ uniformly for actual endpoint foci.
If $a\ge a_0$, then for a sufficiently small fixed $t$ in \eqref{eq:25}, we obtain
$H\ge a_0t/4$. Therefore, before any
target-segment Taylor expansion is made, the smallness of $H$ implies that $0<a<a_0$.
The canonical identities then yield, on the controlled product, 
\[
 g_a(q)=aq_n+O(a^2|q|^2),\qquad
 \tfrac12aq_n\le g_a(q)\le\tfrac32aq_n,
\]
where we have used $|q|^2\le Cq_n$. Consequently, we have $H\le Ca^2$ and
$S_H(x)\subset\{q_n<Ca\}\cap K_{C_0a^2}$.
At a highest point of the set $K_{\kappa a^2}$, we have 
$g_a-v\ge(c\sqrt\kappa/2-\kappa)a^2$.
Choosing $\kappa$ sufficiently small then yields $H\approx a^2$.
At the maximum point $q_x$, \eqref{eq:25} now implies
$q_n(x)\approx a$, and localisation gives $|q_x'|\le Ca\ell(a)$.
Therefore $\rho(q_x')=O(a^2\ell(a)^2)=o(a)$, which proves the distance
comparison and the asserted enclosing boundary strip.

Changing the focus does not destroy the enclosing slide. Set
$\Psi_a(q)=c_\eta(q,ae_n)$ and $B=c_{q\eta}(q_x,ae_n)$. Then
\[
 \begin{split}
 &Q(q)=B^{-T}(\Psi_a(q)-\Psi_a(q_x)),\qquad
 \Psi_a(q)=q+O(a|q|^2),\\
 &D\Psi_a=I+O(a|q|),\qquad B=I+O(a^2\ell(a)).
 \end{split}
\]
For $T=T_{C_0a^2}$ and $q\in S_H(x)$, we have
$|a^{-1}TQ(q)-a^{-1}T(q-q_x)|\le Ca^2\ell(a)^3=o(1)$,
hence $Q(S_H(x))\subset aT^{-1}B_C$.

Let $h=3H/4$, $L_0=aT^{-1}$, $P_0=hL_0^{-T}$, and normalise
$S_h(x)$ with $L_0$ and height $h$. The convex image $D_0$ lies in
$B_C$, and the normalised potential $U_0$ has minimum zero and
boundary value one. The associated normalised double canonical cost
$d_0(\xi,\zeta)=h^{-1}c(L_0\xi,P_0\zeta)$ satisfies
\[
 \|d_0-\xi\cdot\zeta\|_{C^0(B_C\times B_1)}
 \le Ch\|L_0\|^2\|L_0^{-1}\|^2\le Ca^2\ell(a)^4\to0.
\]
Fix a sufficiently small $r$ so that $|d_0|<1/2$ on $D_0\times B_r$.
For every $\zeta\in B_r$, the function $U_0-d_0(\cdot,\zeta)$ has an interior
minimum, since $U_0-d_0(\cdot,\zeta)=0$. The twist condition identifies $\zeta$
as an actual transport value. The Jacobian of this map on the graph is bounded
above by $C(\det L_0)^2/h^n\le C$, so the area formula gives
$|D_0|\ge c$. A convex subset of $B_C$ with this volume has all
its John semi-axes bounded above and below by positive constants. Reserving the normalisation by $L_0$ then yields  the
remaining estimates in \eqref{eq:24}.
\end{proof}

Denote $Q(S_h(x))$ by $b+LD$, where $B_1\subset D\subset B_n$,
and set $P=hL^{-T}$. In these coordinates, define
\[
 d_h(\xi,\zeta)=h^{-1}\{c(b+L\xi,P\zeta)-c(b,P\zeta)\},\qquad
 U_h(\xi)=h^{-1}w_x(Q^{-1}(b+L\xi)).
\]
Here $d_h$ denotes a normalised cost only within this subsection; the normal
height of Theorem~\ref{thm:1} plays no role. The double canonical
cost admits the expansion $c(Q,Z)=Q\cdot Z+R(Q,Z)$, where
$|D_Q^iD_Z^jR|\le C|Q|^{2-i}|Z|^{2-j}$ for $0\le i,j\le2$.
Moreover, $|D_Q^3R|\le C|Z|^2$ and $|D_Z^3R|\le C|Q|^2$.
Since $|b|\le n\|L\|$, the chain rule yields, on every fixed product,
\begin{equation}\label{eq:26}
 \|d_h-\xi\cdot\zeta\|_{C^{2,1}(B_{R_*}\times B_{R_*})}
 \le C_{R_*}h^{-1}\|L\|^2\|P\|^2
 \le C_{R_*}H\ell(H)^4\longrightarrow0.
\end{equation}
For pure third derivatives, the estimate acquires an additional factor
$\|L\|$ or $\|P\|$, each of which is at most one. The phase neighbourhood
ensures that the full products are admissible, including source points outside
$\Omega$. If $J(Q)=|\det D_QX|$, then
\begin{equation}\label{eq:27}
 \begin{split}
 &\det(D^2U_h-A_h(\xi,DU_h))=F_h(\xi),\quad U_h|_{\partial D}=1,
 \quad\min U_h=0,\\
 &F_h=\frac{(\det L)^2}{h^n}J(b+L\xi)^2f(X(b+L\xi)),\quad
 0<\lambda_2\le F_h\le\Lambda_2,\\
 &\operatorname{osc}_DF_h\le
 C\{\omega_f(C\sqrt H\ell(H))+\sqrt H\ell(H)\}\longrightarrow0.
 \end{split}
\end{equation}
The directional inequality and Lemma~\ref{lem:6} remain valid, with
$\kappa_h\le C\|L\|\to0$.

\begin{lemma}\label{lem:11}
For every fixed $0<\sigma<1$ and $p>1$, all sufficiently small $H$ satisfy
\begin{equation}\label{eq:28}
 \int_{\{U_h<\sigma\}}|D^2U_h|^p\le C_{p,\sigma},\qquad
 B_{b_\sigma}(\xi_h)\subset\{U_h<\sigma\},\quad U_h(\xi_h)=0.
\end{equation}
\end{lemma}
\begin{proof}
Take any sequence $H_j\to0$. Convex-envelope compactness and testing
the original equations as in Lemma~\ref{lem:7} give a convex limit
$V$ on $B_1\subset D\subset B_n$ satisfying $\det D^2V=F_0$,
where $F_0$ is a positive constant. To identify the trace, choose a smooth
uniformly convex outer domain $D^\delta$ together with the classical solution of
$\det D^2V^\delta=2\Lambda_2$ with $V^\delta=1$ on $\partial D^{\delta}$.
For fixed $\delta$, $V^{\delta}$ is a strict augmented subsolution for large $j$
by \eqref{eq:26}, so that $U_j\ge V^\delta$. The uniform Alexandrov
estimate $V^\delta\ge1-C\operatorname{dist}(\cdot,\partial D^\delta)^{1/n}$
implies that $V=1$ continuously on $\partial D$. This estimate also ensures that every smaller
sublevel set remains away from $\partial D_j$, including the minima. Consequently,  $\min V=0$.

On smaller sublevels, $V$ is uniformly smooth and strictly convex.
Indeed, convexity bounds the gradient of $V$ away from the boundary. We apply
the classical Pogorelov estimate to $V-\sigma_3$ on the intermediate
section $\{V<\sigma_3\}$, and then restrict to $\{V\le\sigma_2\}$,
where $\sigma_2<\sigma_3<1$. The determinant equation supplies the lower
Hessian bound, and classical interior estimates yield uniform higher-order estimates. The justification of these estimates  for the
initial Alexandrov solution follows from smooth outer approximation. Lemma~\ref{lem:6} implies that $U_j$ is uniformly
semiconvex, with the negative Hessian bound vanishing on compact subsets.
Difference quotients therefore yield the convergence $DU_j\to DV$ locally uniformly.
The actual maps $T_j$ converge there to $DV$: on a fixed sufficiently
large target ball, the twist map $(d_j)_\xi$ is uniformly close to the identity,
and its image contains the bounded gradients. A local inverse is obtained
by Brouwer's theorem; the global twist condition then identifies this local inverse with the actual transport value.

Fix a center $\xi_0$ in a smaller sublevel and choose uniform $r,m,M>0$ such that
$B_{4r}(\xi_0)$ lies inside a larger sublevel and $mI\le D^2V\le MI$ there.
For $t_j=T_j(\xi_0)$ and a small fixed $\rho<mr/8$, the function
$U_j-d_j(\cdot,y)$ attains its minimum over $\overline B_{2r}(\xi_0)$
at an interior point whenever $y\in B_\rho(t_j)$. This follows from the
quadratic lower separation of $V$ on $\partial B_{2r}(\xi_0)$ and local uniform
convergence established above. At such a minimum, the twist condition yields $T_j(\xi)=y$, and the injectivity
of  Lemma~\ref{lem:1} excludes any other preimage. Thus, using the convergence $T_j\to DV$ also for the inner inclusion, the actual
transport can be restricted to
\[
 Y_j=B_\rho(t_j),\qquad X_j=T_j^{-1}(Y_j),\qquad
 B_{r_0}(\xi_0)\subset X_j\subset B_{2r}(\xi_0)
\]
for a fixed $r_0>0$.

Set $s_j=F_j(\xi_0)^{1/n}$, and let $\xi=\xi_0+q$, $y=t_j+s_j\zeta$.
We define the cost and potential by
\[
 \begin{aligned}
 s_j\widehat d_j(q,\zeta)
    &=d_j(\xi,y)-d_j(\xi,t_j)-d_j(\xi_0,y)+d_j(\xi_0,t_j),\\
 s_j\widehat U_j(q)
    &=U_j(\xi)-U_j(\xi_0)-d_j(\xi,t_j)+d_j(\xi_0,t_j).
 \end{aligned}
\]
The four cost terms combine to yield $s_j q\cdot\zeta$ for $d_j(\xi,y)=\xi\cdot y$.
Since $s_j$ is bounded above and away from zero, \eqref{eq:26} gives
$\widehat d_j\to q\cdot\zeta$ in $C^{2,1}$ on each fixed product.
The new domains lie between $B_{1/K}$ and $B_K$, and the target domain is a ball.
The associated map is given by $\widehat T_j(q)=(T_j(\xi)-t_j)/s_j$, with densities
\begin{equation}\label{eq:29}
 \widehat\rho_j(q)=
 \frac{F_j(\xi_0+q)}{F_j(\xi_0)|\det(\widehat d_j)_{q\zeta}
       (q,\widehat T_j(q))|},\qquad \widehat\rho_j^*=1.
\end{equation}
Indeed, division by $s_j$ scales the augmented determinant by $s_j^n$,
while the target dilation cancels in the mixed cost derivative.
Injectivity and the area formula ensure the equality of total masses and the
transport identity. Equations~\eqref{eq:26}--\eqref{eq:27} imply
$\|\widehat\rho_j-1\|_\infty\to0$. Every source point retains its
original support in the restricted target; moreover, the supremum of these support allows the potential to be extended if an ambient extension is required.
All hypotheses of \cite[Proposition 2.4]{ChenFigalli2017}  are now satisfied
for the minimizing cost $-\widehat d_j$. Its proof in Section 3.4
yields a uniform local $L^p$ bound for $D^2\widehat U_j$.
Since $D^2U_j(\xi)=s_jD^2\widehat U_j(q)+(d_j)_{\xi\xi}(\xi,t_j)$,
\eqref{eq:26} gives the same bound for $D^2U_j$. If the total mass is normalised to one, a common
bounded dilation of source and target, together with division of cost
and potential by the square of the dilation factor, reduces the problem to this convention.
This normalisation preserves the density one on the target.
A finite covering argument yields the first assertion in \eqref{eq:28}.
The minima converge to the unique minimum of $V$; its uniform upper
Hessian bound and the local uniform convergence yield the second assertion.
All choices are uniform, since a sequence violating them would admit
the compactness limit described above, which is a contradiction.
\end{proof}

\begin{proof}[Proof of Theorem~\ref{thm:2}]
We apply Lemma~\ref{lem:11} with $h=3H/4$ and $\sigma=2/3$.
On this smaller section, the normalised gradient and Hessian matrix are bounded,
so that \eqref{eq:28} also provides an $L^p$ bound for $W_U=D^2U-A_h(\xi,DU)$.
The associated covariance matrix is given by $W_Q=hL^{-T}W_UL^{-1}$.
The canonical coordinate Jacobians are bounded in both directions. Consequently,
\eqref{eq:24} yields
\begin{equation}\label{eq:30}
 B_{b\sqrt H/\ell(H)}(x)\subset S_{H/2}(x),\qquad
 \int_{S_{H/2}(x)}|D^2u|^p\le C_pH^{n/2}\ell(H)^{2p}.
\end{equation}
Here we have used the estimate 
\[
 h^p\|L^{-1}\|^{2p}|\det L|\le CH^{n/2}\ell(H)^{2p}
\]
together with the uniform bound for $A(x,Du)$. The threshold for $H$
is uniform over all centers, as a consequence of  \eqref{eq:26}--\eqref{eq:27}.

For $E_r=\{r<\operatorname{dist}(x,\partial\Omega)<2r\}$, set
$\varrho_r=br/\ell(r)$. Choose a maximal disjoint family of balls
$B_{\varrho_r/4}(x_i)$ with centers in $E_r$.
These balls lie in a boundary strip of volume at most $Cr$, so that their
number is bounded by $Cr^{1-n}\ell(r)^n$. The doubled balls $B_{\rho_r/2}(x_i)$ cover $E_r$
and are contained in the corresponding sections from \eqref{eq:30}, since
$H(x_i)\approx r^2$. Hence
\[
 \int_{E_r}|D^2u|^p\le C_pr\ell(r)^{n+2p}.
\]
Summing over the dyadic scales $r=2^{-k}$, we obtain a convergent series. Proposition~\ref{prop:1} handles the remaining
fixed interior rigion, while the initial $C^1$ bound controls the lower-order
terms. This completes the proof of the theorem.
\end{proof}

\section{The global H\"older estimate}

For $f\in C^\alpha(\overline\Omega)$, we first establish a quadratic
expansion at each boundary point. In canonical coordinates, the matrix
term vanishes to second order in the gradient, and Theorem~\ref{thm:1}
controls the gradient near the point of contact. A comparison with the
constant-density solutions of the Monge-Amp\`ere equation then yields,
by iteration, a remainder of order $2+\alpha$. To pass from the
boundary expansions to the global H\"older estimate, we combine
Theorem~\ref{thm:2} with Morrey's inequality, which controls the variation
of the matrix on interior balls.

\subsection{The matrix along the solution}
At a boundary point $z$, we retain the canonical coordinates $q=Q_z(x)$, the
support $m_z$, and the remainder $v_z$. The inner localisation inclusion
and the logarithmic slide bound then imply
\begin{equation}\label{eq:31}
 0\le v_z(q)\le C|q|^2\ell(|q|)^2.
\end{equation}
Indeed, for $|q|=r$, we take $h=C_1r^2\ell(r)^2$ so that
$|T_hq|<c\sqrt h$. The same growth estimate holds in physical
coordinates, since the coordinate  maps are uniformly bi-Lipschitz.

\begin{lemma}\label{lem:12}
For a neighborhood of each boundary point, we have the following formulas holding uniformly:
\begin{equation}\label{eq:32}
 \begin{split}
 &|Du(x)-Du(z)|\le C|x-z|\ell(|x-z|)^2,\\
 &|Dv_z(q)|\le C|q|\ell(|q|)^2,\qquad
 |a_z(q)|\le C|q|^2\ell(|q|)^4,
 \quad a_z(q)=\widehat A_z(q,Dv_z(q)).
 \end{split}
\end{equation}
\end{lemma}
\begin{proof}
For boundary points $z,w$ at distance $r$, let
$l_z(x)=u(z)+Du(z)\cdot(x-z)$ and define $l_w$ analogously.
By \eqref{eq:31} and the uniform Hessian bound for the supports,
we have $|u-l_z|+|u-l_w|\le Cr^2\ell(r)^2$ on a common interior ball of
radius $cr$. The gradient of the affine function $l_z-l_w$ is therefore bounded by
$Cr\ell(r)^2$. This proves the boundary-pair assertion.
For an interior point $x$ at distance $d$ from a nearest boundary point
$z_x$, the estimate \eqref{eq:31} bounds the oscillation of $v_{z_x}$ on the
canonical image of $B_{d/2}(x)$ by $Cd^2\ell(d)^2$. This image
contains a ball of radius $cd$ centered at $Q_{z_x}(x)$.
The two-sided directional estimate \eqref{eq:3} then bounds the gradient
there by the oscillation divided by $d$. Returning to physical coordinates
yields $|Du(x)-Du(z_x)|\le Cd\ell(d)^2$.
Combining this with the boundary-pair estimate and the inequality $|z_x-z|\le2|x-z|$, we obtain the first line.
The chain rule yields the second line, in view of  the canonical quadratic
vanishing of $\widehat A_z$ at zero gradient.
\end{proof}

Set $L_h=\sqrt h T_h^{-1}$ and $U_h=h^{-1}v_z(L_h\xi)$, and regard
\[
 a_h(\xi)=h^{-1}L_h^Ta_z(L_h\xi)L_h
\]
as a prescribed matrix field. Equation~\eqref{eq:32} implies that on fixed
normalised neighbourhoods,
\begin{equation}\label{eq:33}
 |a_h(\xi)|\le Ch\ell(h)^8|\xi|^2\ell(|\xi|)^4,
 \qquad \sup_{0<|\xi|\le1}\frac{|a_h(\xi)|}{|\xi|^\alpha}\to0.
\end{equation}
Here, the inequality $\ell(|L_h\xi|)\le C\ell(h)\ell(|\xi|)$ follows from the upper
and lower singular-value bounds for $L_h$.
Let $F_h(\xi)=\widetilde f(L_h\xi)$ denote the density of $U_h$, and set
$s_h=F_h(0)^{-1/(2n)}$. Define $\bar U_h(\xi)=U_h(s_h\xi)$ on the
dilated section. The density of $\bar U_h$ is given by $f_h(\xi)=s_h^{2n}F_h(s_h\xi)$, so that
$f_h(0)=1$ and $|f_h(\xi)-1|\le o(1)|\xi|^\alpha$.
Both $s_h$ and $s_h^{-1}$ are uniformly bounded. This dilation
preserves the unit section height and the estimate \eqref{eq:33}, with the matrix replaced by
$s_h^2a_h(s_h\xi)$.
The transformed boundary and its trace admit quadratic Taylor expansions
with errors $o(1)|\xi'|^{2+\alpha}$; indeed, their third-order errors carry a
factor $\sqrt h$ times a fixed power of $\ell(h)$, which tends to zero.
To be noticed, all these smallness statements are uniform.

\begin{lemma}\label{lem:13}
Let $D\subset B_R(x_*)$ and let $U\in C^2(D)\cap C(\overline D)$ satisfy
\[
 W=D^2U-a>0,\qquad \det W=f>0\quad\hbox{in }D.
\]
Suppose that $V$ is a continuous convex
Alexandrov solution of $\det D^2V=b_0>0$ on $D$, and assume that
\[
 \begin{gathered}
  \|a\|\le\eta,\qquad |f^{1/n}-b_0^{1/n}|\le\epsilon
       \quad\hbox{in }D,\\
  |U-V|\le\delta \quad\hbox{on }\partial D.
 \end{gathered}
\]
Then
\begin{equation}\label{eq:34}
 \|U-V\|_{L^\infty(D)}\le\delta+\tfrac12R^2(\eta+\epsilon).
\end{equation}
\end{lemma}

\begin{proof}
Fix $\tau>0$ and set
\[
 Q(x)=|x-x_*|^2-R^2,\qquad t=\tfrac12(\eta+\epsilon+\tau).
\]
We claim that
\[
 U+tQ-\delta\le V\le U-tQ+\delta\quad\hbox{in }D.
\]
Both inequalities hold on $\partial D$, since $Q\le0$ and
$|U-V|\le\delta$ on $\partial D$. To establish them in the interior, we use the convex-viscosity characterization
of the Alexandrov equation satisfied by $V$.

Suppose, for contradiction, that the first inequality fails, then $U+tQ-V-\delta$ attains a positive
maximum at some interior point, and a vertical translate of $U+tQ$ therefore
touches $V$ from below at that point. Since $-\eta I\le a\le\eta I$,
the Hessian of $U+tQ$ satisfies
\[
 D^2U+2tI=W+a+2tI\ge W+(\epsilon+\tau)I>0,
\]
 where we have used the definition of $t$. The viscosity supersolution inequality combined with Minkowski's determinant
inequality then yields, at the contact point,
\[
 b_0^{1/n}\ge\bigl[\det(D^2U+2tI)\bigr]^{1/n}
 \ge f^{1/n}+\epsilon+\tau>b_0^{1/n},
\]
where the strict inequality follows from the assumption
$|f^{1/n}-b_0^{1/n}|\le\epsilon$. This is a contradiction.

Suppose now that the second inequality fails, then $U-tQ-V+\delta$ attains a negative
minimum at some interior point, and a vertical translate of $U-tQ$ therefore touches
$V$ from above at that point. The convexity of $V$ implies $D^2U-2tI\ge0$ at
the contact point, which follows by considering second differences in each
direction. The viscosity subsolution inequality for $V$ gives
$\det(D^2U-2tI)\ge b_0$. On the other hand,
\[
 W=(D^2U-2tI)+(2tI-a)\ge(D^2U-2tI)+(\epsilon+\tau)I.
\]
Applying the Minkowski's determinant inequality once again, we obtain
\[
 f^{1/n}\ge\bigl[\det(D^2U-2tI)\bigr]^{1/n}+\epsilon+\tau
 \ge b_0^{1/n}+\epsilon+\tau>f^{1/n},
\]
which is again a contradiction. Hence the claimed comparison holds, and for all $x\in D$
\[
 |U(x)-V(x)|\le\delta-tQ(x)
 \le\delta+\tfrac12R^2(\eta+\epsilon+\tau)\quad(x\in D).
\]
Let $\tau$ decreases monotonically to 0, we obtain \eqref{eq:34}.
\end{proof}

\subsection{Boundary quadratic approximation}
We first establish a uniform starting scale. We apply the compactness proof of
Lemma~\ref{lem:7}  to $M_h=s_h\sqrt h T_h^{-1}$, since
localisation provides uniform inner and outer balls, we have $\det M_h=s_h^nh^{n/2}$,
and $h\|M_h\|^2\|M_h^{-1}\|^2\le Ch\ell(h)^4\to0$.
For $h_j$ decreasing monotonically to 0, set $U_j=\bar U_{h_j}$ and let $\Gamma_j$
denote the lower convex envelope of $U_j$. Along a subsequence, a limit $V$ on
the half-domain $O$ satisfies $\det D^2V=1$, since $f_{h_j}\to1$
uniformly on bounded sets. The normalised physical boundaries converge to
$\{\xi_n=0\}$, and their traces are given by
$p_j(\xi')+o(1)$, where $p_j$ are homogeneous quadratics satisfying
$c|\xi'|^2\le p_j\le C|\xi'|^2$. Passing to a further subsequence, we may assume that $p_j\to p_\infty$.
The limiting closed trace domain is exactly
$G=\{p_\infty\le1\}\cap\{\xi_n=0\}$. The boundary datum is then
$p_\infty$ on $G$ and 1 on the remainder of $\partial O$.
These are  precisely the hypotheses of
\cite[Theorem 6.4]{Savin2013Boundary}, with the closed-epigraph convention.
The theorem yields a uniform $C^{3,1}$ bound on a smaller half-ball.

We also require convergence up to the flat part of the boundary. On a disk of uniform small radius, the affine tangent plane to $p_\infty$ at $\xi'$ lies below
$p_\infty$ on $G$ and below the value one on the remainder of the bounded boundary.
Thus the convexification of the trace equals $p_\infty$ on that disk.
For points $\xi_j$ converging to the disk, we join the physical boundary
point lying directly below $\xi_j$ to a point at a fixed small height above it.
The convexity of $\Gamma_j$  then bounds $\Gamma_j(\xi_j)$ by the corresponding interpolation
of these endpoint values. Moreover, the coefficient of the upper endpoint tends
to zero as $\xi_j$ approaches the disk. The trace expansion, together with $\|U_j-\Gamma_j\|_\infty\to0$
yields the required upper limit, while the closed-epigraph convergence provides the
lower limit. Together with the interior convergence, this yields uniform
convergence on a small moving closed half-neighbourhood.

The Taylor polynomial of $V$ at the origin is the sum of a linear function  normal to the boundary and a
positive definite quadratic form.  We perform a fixed further small dilation, so
that the cubic error is below a chosen quadratic tolerance. The uniform convergence transfers this approximation to $U_j$.
We then subtract the normal linear part, adjust the tangential quadratic form to match
the current boundary trace, and finally adjust the normal-normal coefficient 
to restore the determinant to one. The positive definiteness of the tangential block ensures that this
last adjustment is both unique and small. In view of \eqref{eq:33}, for
any fixed small $\delta,\epsilon_0$, we obtain that at a uniform starting
scale, a function $U$ satisfies 
\begin{equation}\label{eq:35}
 \begin{split}
 &D^2U-a>0,\quad\det(D^2U-a)=f,\quad
 |a(x)|+|f(x)-1|\le\delta\epsilon_0|x|^\alpha,\\
 &|U-P_0|\le\epsilon_0\quad\hbox{in }\Omega_0\cap B_1,\qquad
 \det D^2P_0=1,\quad c_0I\le D^2P_0\le C_0I.
 \end{split}
\end{equation}
The convex domain lies above its tangent plane, with boundary given by
$x_n=g(x')\ge0$. The homogeneous quadratics $q$ and $p$ then satisfy
\begin{equation}\label{eq:36}
 \begin{split}
 &|g(x')-q(x')|\le\delta\epsilon_0|x'|^{2+\alpha},\quad
 \|q\|\le\delta\epsilon_0,\\
 &|U(x',g(x'))-p(x')|\le\delta\epsilon_0|x'|^{2+\alpha},\quad
 p(x')=P_0(x',0),\quad c_0|x'|^2\le p(x')\le C_0|x'|^2.
 \end{split}
\end{equation}
Subtracting linear functions leaves the matrix field unchanged along the original solution. Moreover, at this fixed starting scale, all
coordinate maps, their inverses, and the scalar factors have uniform bounds.

\begin{proposition}\label{prop:3}
For sufficiently small $\delta,\epsilon_0$, depending only on
$n,\alpha,c_0,C_0$, \eqref{eq:35}--\eqref{eq:36} imply
\begin{equation}\label{eq:37}
 |U(x)-l(x)-P(x)|\le C\epsilon_0|x|^{2+\alpha}
 \quad(x\in\Omega_0\cap B_c),
\end{equation}
where $l$ is normal linear and $P$ is quadratic with uniformly
positive definite Hessian and bounded coefficients.
\end{proposition}
\begin{proof}
Fix $0<\tau<1$. At the scale $r=\tau^m$, we construct $l_m=\gamma_mx_n$ and a
homogeneous quadratic $P_m$ such that
\[
 \begin{split}
 &|U-l_m-P_m|\le\epsilon_0r^{2+\alpha}\quad\hbox{on }\Omega_0\cap B_r,\\
 &P_m(x',0)=p(x')-\gamma_mq(x'),\qquad
 \det D^2P_m=1,\quad cI\le D^2P_m\le CI,\quad|\gamma_m|\le1.
 \end{split}
\]
These conditions hold initially.  We then set $v(X)=r^{-2}(U-l_m)(rX)$,
$\epsilon=\epsilon_0r^\alpha$, and $\Omega_r=r^{-1}\Omega_0$.
,then we have $|v-P_m|\le\epsilon$, the matrix and density errors are bounded by
$\delta\epsilon$, and on the physical boundary within $B_1$, we have
\begin{equation}\label{eq:38}
 0\le X_n\le2\delta\epsilon,\qquad
 |v-P_m|\le C_1\delta\epsilon.
\end{equation}
For the latter inequality, the tangential identity eliminates the term
$\gamma_mq$ while the remaining graph errors and the mixed terms in $P_m$
are bounded by using $r\le r^\alpha$.

We first control the layer near the boundary. On the annulus
$1/4<|X|<1/2$, we choose
$\phi(X)=c_\beta(4^\beta-|X|^{-\beta})$, normalised so that $\phi=1$ at
$|X|=1/2$. For $\beta$ fixed large, the contraction of $D^2\phi$ with
every matrix in the fixed ellipticity class is strictly negative.
For $Y=(Y',-1/4)$, we use
$P_m\pm\epsilon(C_1\delta+\phi(X-Y))$ as barriers on the
intersection of this annulus with $\Omega_r$. For $\epsilon_0$ small, the Hessians of these barriers are positive, and their
determinants lie below $1-c\epsilon$ and above
$1+c\epsilon$, respectively. Subtracting the actual field changes these determinants
by at most $C\delta\epsilon$. The boundary comparison follows from
\eqref{eq:38} together with $|v-P_m|\le\epsilon$ on the outer sphere, since the
inner sphere has no points of $\Omega_r$ except for a possible tangency.
The maximum principle therefore applies on the elliptic branch,
and varying $Y'$ yoelds
\begin{equation}\label{eq:39}
 |v(X)-P_m(X)|\le\epsilon(C\delta+CX_n)
 \quad\hbox{in }B_{1/8}\cap\Omega_r.
\end{equation}

Fix a fixed smaller half-ball $D=B_{r_*}^+$, $r_*<1/8$, and
$t_0=3\delta\epsilon$. On its hemisphere, we prescribe $P_m$ below
$t_0$, $v$ above $2t_0$, and interpolate between the two using a continuous cutoff. This strip lies in $\Omega_r$, and \eqref{eq:39} bounds
the corresponding interpolation error by $C\delta\epsilon$. On the flat disk, we
prescribe $P_m$. These continuous data admit their own boundary convex
envelope: hemisphere points are strictly exposed, while any convex
combination yielding a point of the flat disk involves only that disk, where
$P_m$ is convex. Consequently, the convex Alexandrov solution $w$ of
$\det D^2w=1$ attains these data.
The radial barriers constructed above, once multiplied by a fixed sufficiently large
constant, also bound $w$ near the flat disk. On the hemisphere, their
increments dominate $\epsilon(C\delta+CX_n)$ from \eqref{eq:39}, while
on the disk the data coincide exactly with $P_m$. Moreover, at a point lying directly above
$Y'$, their increments are $O(\epsilon(\delta+X_n))$.
In particular, we obtain $|w-P_m|\le C\delta\epsilon$ on the physical graph.
Applying Lemma~\ref{lem:13} on $D\cap\Omega_r$, we now obtain
\begin{equation}\label{eq:40}
 \|v-w\|_\infty\le C\delta\epsilon,\qquad
 \|w-P_m\|_{L^\infty(D)}\le2\epsilon.
\end{equation}
The second bound can also be obtained directly from a comparison with
$P_m\pm2\epsilon$.

Choose a fixed $t>0$ with $4t<\min_{|X|=r_*}P_m(X)$ uniformly in
$m$, and we take $\epsilon_0$ such that $2\epsilon_0<t/2$. The sublevel
$D_t=\{w<t\}\cap D$ stays away from the curved hemisphere and is
contained between half-balls of radii $c\sqrt t$ and $C\sqrt t$ by \eqref{eq:40}.
The function $t^{-1}w(\sqrt t\,X)$ still has
determinant one on $X=t^{-1/2}D_t$. The flat trace is $P_m(X',0)$ on
$\{P_m(X',0)\le1\}$, then the remaining boundary
value is one by homogeneity. Applying \cite[Theorem 6.4]{Savin2013Boundary} and
scaling back, we obtain a uniform $C^{3,1}$ bound on a smaller half-ball. In this region, $w-P_m$ solves a uniformly elliptic
linear equation with coefficients controlled in terms of the ellipticity constants, and with zero boundary
data on the flat part. The linear boundary estimates, together with the $L^\infty$ bound in
\eqref{eq:40}, yield
\begin{equation}\label{eq:41}
 \begin{split}
 &|w-P_m-\widetilde\gamma_mX_n-\widetilde P_m(X)|
       \le C\epsilon|X|^3,\\
 &|\widetilde\gamma_m|+\|D^2\widetilde P_m\|\le C\epsilon,\quad
 \widetilde P_m(X',0)=0,\quad
 \det D^2(P_m+\widetilde P_m)=1.
 \end{split}
\end{equation}
The last equality is precisely the equation satisfied at the origin by the Taylor Hessian.

Set $\gamma_{m+1}=\gamma_m+r\widetilde\gamma_m$ and
\[
 P_{m+1}=P_m+\widetilde P_m-r\widetilde\gamma_mq(x')+s_mx_n^2,
\]
where $s_m$ restores the determinant to one. The positive definiteness of the tangential block
gives the uniqueness of $s_m$ with $|s_m|\le C\delta\epsilon$.
The coefficient increments are bounded by $C\epsilon$ and the slope
increment is bounded by $Cr\epsilon$; moreover, the tangential identity is preserved.
Equations~\eqref{eq:40}--\eqref{eq:41} yield an error of at most
$Cr^2\epsilon(\delta+\tau^3)$ on $B_{r\tau}$.
We first choose $\tau$ so that $C\tau^3\le\tau^{2+\alpha}/4$, then $\delta$ such that $C\delta\le\tau^{2+\alpha}/4$, and finally $\epsilon_0$ so
that the summable increments preserve both ellipticity and the slope bound.
This closes the induction argument. Summing the tails of the coefficient series
$\sum_{j\ge m}C\epsilon_0\tau^{j\alpha}$ and
$\sum_{j\ge m}C\epsilon_0\tau^{j(1+\alpha)}$, we obtain \eqref{eq:37}.
\end{proof}

Undoing the fixed normalisation and the smooth cost coordinates, we obtain,
at every $z\in\partial\Omega$, a quadratic polynomial $R_z$ such that
\begin{equation}\label{eq:42}
 |u(x)-R_z(x)|\le C|x-z|^{2+\alpha},\qquad
 \|D^2R_z\|\le C\quad(x\in\overline\Omega\cap B_c(z)).
\end{equation}
Since $u\in C^2(\overline\Omega)$, the polynomial $R_z$ coincides with the Taylor polynomial of $u$ at $z$.
In canonical coordinates, the linear part of $R_z$ vanishes and its Hessian is
uniformly positive definite. Comparing $R_z-R_w$ on a common interior
ball of radius comparable to $|z-w|$, and using  the equivalence of norms
for quadratic polynomials, we obtain
\begin{equation}\label{eq:43}
 |D^2R_z-D^2R_w|\le C|z-w|^\alpha.
\end{equation}

\subsection{Interior approximation and the global bound}
\begin{lemma}\label{lem:14}
Suppose that $w\in C^2(B_1)$, and that $b$ is a continuous symmetric matrix field satisfying
$D^2w-b>0$, $\det(D^2w-b)=g$, $b(0)=0$, $g(0)=1$.
Let $P(X)=\tfrac12X^TBX$, where $\det B=1$ and $\lambda_*I\le B\le\Lambda_*I$.
There exist constants $\delta,E_*,C>0$, depending only on $n,\alpha,\lambda_*,\Lambda_*$,
such that
\[
 \|w-P\|_\infty\le E\le E_*,\qquad
 [b]_{C^\alpha(B_1)}+[g]_{C^\alpha(B_1)}\le\delta E
\]
imply that
\begin{equation}\label{eq:44}
 \|D^2w-B\|_{L^\infty(B_{1/2})}+[D^2w]_{C^\alpha(B_{1/2})}\le CE.
\end{equation}
\end{lemma}
\begin{proof}
We first construct a Taylor expansion at the origin. At the scale $r=\tau^m$, assume that
$|w-l_m-P_m|\le Er^{2+\alpha}$ on $B_r$, where $l_m$ is affine,
$\det D^2P_m=1$, and $D^2P_m$ lies in a fixed ellipticity range.
Setting $v=r^{-2}(w-l_m)(rX)$, and we let $\epsilon=Er^\alpha$.
Then $|v-P_m|\le\epsilon$, and both the prescribed field and the density
errors are bounded by $\delta\epsilon|X|^\alpha$.
We solve $\det D^2V=1$ on $B_{1/2}$ with boundary data $v$.
Since the sphere is strictly convex, arbitrary continuous boundary data are admissible.
Classical comparison and Lemma~\ref{lem:13} give
$\|V-P_m\|_\infty\le\epsilon$ and
$\|v-V\|_\infty\le C\delta\epsilon$.

We recall the uniform model estimates, whereby convexity bounds $DV$ on
$B_{3/8}$. Two-sided difference quotients then give
$|p-DP_m(X)|\le C\sqrt\epsilon$ for every $p\in\partial V(X)$, on a smaller ball. The closeness to $P_m$ therefore provides, at every point
of $B_{1/4}$, a supporting section of fixed positive height, contained
in $B_{3/8}$ and containing a fixed ball. The classical Pogorelov
estimate applied to these sections bounds $D^2V$ from above, while the determinant equation provides the lower bound. Smooth approximation and the classical interior
estimates give uniform $C^{3,\gamma}$ bounds on a smaller ball.
The difference $Z=V-P_m$ then solves
\[
 \left(\int_0^1\operatorname{cof}((1-t)D^2P_m+tD^2V)\,dt\right):D^2Z=0.
\]
The linear estimates preserve the bound $\|Z\|_\infty\le\epsilon$, so that the coefficients in the quadratic Taylor expansion of $Z$ are of size $C\epsilon$ with remainder $C\epsilon|X|^3$. The Hessian of the resulting quadratic, added to $D^2P_m$, has determinant one.
We add the constant and linear terms, scaled by $r^2$ and $r$ respectively, to $l_m$,
and add the corresponding homogeneous quadratic to $P_m$. The resulting error at the scale $r\tau$ is
at most $Cr^2\epsilon(\delta+\tau^3)$. Choosing $\tau,\delta,E_*$
as in Proposition~\ref{prop:3} closes the induction and yields an
expansion with remainder $CE|X|^{2+\alpha}$, whose Hessian lies within $CE$ of $B$.

For each $Y\in B_{1/2}$, we recenter using the fixed radius $\rho=1/4$.
We subtract the quadratic $\tfrac12X^Tb(Y)X$, divide by $g(Y)^{1/n}$,
and adjust the comparison quadratic by a scalar multiple of $|X|^2$
to restore the determinant to one. Since $|b(Y)|+|g(Y)-1|\le C\delta E$,
these adjustments are of order $O(\delta E)$, and they preserve the ellipticity.
The new flatness is at most $C\rho^{-2}E$ while the coefficient seminorms
at most $C\rho^\alpha\delta E$. Reducing the original thresholds
by fixed factors makes the preceding iteration uniform in $Y$.
Consequently, every point has its own Taylor polynomial $Q_Y$, with remainder
$CE|X-Y|^{2+\alpha}$ on a fixed ball and with Hessian within $CE$ of $B$.
For nearby points $Y,Z$, we compare $Q_Y-Q_Z$ on a ball of radius comparable
to $|Y-Z|$ contained in the overlap. Polynomial norm equivalence gives
$|D^2Q_Y-D^2Q_Z|\le CE|Y-Z|^\alpha$.
For separated points, we instead use the bound relative to $B$. This establishes
\eqref{eq:44}, including the factor $E$.
\end{proof}

\begin{proof}[Proof of Theorem~\ref{thm:3}]
We choose $\alpha<\beta<1$ and $p>n/(1-\beta)$. Theorem~\ref{thm:2}
and Morrey's inequality yield
\begin{equation}\label{eq:45}
 \|Du\|_{C^{0,\beta}(\overline\Omega)}\le C.
\end{equation}
Let $x_0$ be a point at small boundary distance $d$, and let $z$ denote its nearest boundary point. In the canonical coordinates at $z$, we set $q_0=Q_z(x_0)$ and $r=\kappa d$,
where $\kappa$ is uniform and small enough that $q_0+rB_4$ is contained in the interior.
We denote $v=v_z$, $a(q)=\widehat A_z(q,Dv(q))$, and let
$Q(q)=\tfrac12q^TH_zq$ be the canonical boundary Taylor polynomial.
We denote the canonical density by
$f_z(q)=(\det DX_z(q))^2f(X_z(q))$, where $X_z=Q_z^{-1}$.
Recall that the Hessian $H_z$ is uniformly positive definite. Equation~\eqref{eq:42} gives
\begin{equation}\label{eq:46}
 \left\|\frac{v(q_0+rX)-Q(q_0)-rDQ(q_0)\cdot X}{r^2}
       -\tfrac12X^TH_zX\right\|_{L^\infty(B_2)}\le Cd^\alpha.
\end{equation}
Moreover, $0\le v\le Cd^2$ on $q_0+rB_4$. The uniform semiconvexity
and two-sided difference quotients of length comparable to $d$ give
$|Dv|\le Cd$ on $q_0+rB_2$. Recall that the canonical identities imply
$|\widehat A_p(q,p)|\le C|p|$ and $|\widehat A_q(q,p)|\le C|p|^2$.
In view of \eqref{eq:45}, these together yield
\begin{equation}\label{eq:47}
 |a|\le Cd^2,\qquad
 [a(q_0+r\,\cdot)]_{C^{0,\beta}(B_2)}\le C(d^{1+\beta}+d^3).
\end{equation}
The scaled density has $C^\alpha$ seminorm at most $Cd^\alpha$.

We freeze $a_0=a(q_0)$ and subtract $\tfrac12X^Ta_0X$ from the
scaled function in \eqref{eq:46}. Its residual matrix is then
$a(q_0+rX)-a_0$. Since the boundary equation yields
$\det H_z=f_z(0)$, we have
$|\det(H_z-a_0)-f_z(q_0)|\le Cd^\alpha$.
An adjustment $tI$ with $|t|\le Cd^\alpha$ makes the comparison
quadratic have determinant $f_z(q_0)$. We then divide both the function and the matrix by
$f_z(q_0)^{1/n}$. Setting $E=Kd^\alpha$, we choose $K$ fixed and large
enough that the flatness is at most $E$ and the density seminorm
is at most $\delta E/2$. By \eqref{eq:47}, the residual matrix seminorm
is also at most $\delta E/2$, provided the collar is sufficiently thin and uniform.
Note that $1+\beta>\alpha$, and the $C^\beta$ bound therefore controls the $C^\alpha$ seminorm
on the fixed ball. Lemma~\ref{lem:14} therefore yields, after undoing
the scalar and quadratic corrections,
\begin{equation}\label{eq:48}
 |D^2v(q_0)-H_z|\le Cd^\alpha,\qquad
 [D^2v]_{C^\alpha(B_{r/2}(q_0))}\le C.
\end{equation}
The factor $r^{-\alpha}$ in the seminorm exactly cancels the factor $E=Kd^\alpha$.

Returning to physical coordinates, we use the exact formula
\[
 D_x^2u=D_x^2m_z+(DQ_z)^TD_q^2v_zDQ_z
          +\sum_a(v_z)_{q_a}D_x^2(Q_z)_a.
\]
The coordinate coefficients are uniformly smooth, $|Dv_z(q_0)|\le Cd$,
and \eqref{eq:45} controls the $C^\alpha$ seminorm of $Dv_z$. Consequently,
\begin{equation}\label{eq:49}
 |D^2u(x_0)-D^2R_z|\le Cd^\alpha,\qquad
 [D^2u]_{C^\alpha(B_{cd}(x_0))}\le C.
\end{equation}
For two collar points separated by less than a fixed fraction of the
smaller boundary distance, we use this local estimate. Otherwise, both
boundary distances are bounded by a constant times their separation, and we
compare each Hessian with its nearest boundary polynomial using
\eqref{eq:43}. This proves the full estimate in the collar.

On a fixed interior compact set, we use Proposition~\ref{prop:1}.
In the double canonical coordinates at an interior contact point, we set
$U_h(X)=h^{-1}v(L_hX)$ and $P_h=hL_h^{-T}$.
Their limits on fixed balls are uniformly positive definite quadratic forms.
The directional estimate on a larger ball yields a bound on  $DU_h$, so that
\[
 \sup_{B_2}|Dv(L_hX)|\le C\|P_h\|,\qquad
 [Dv(L_h\,\cdot)]_{C^\beta(B_2)}\le C\|L_h\|^\beta.
\]
For the actual matrix and density
\[
 a_h=h^{-1}L_h^TA(L_hX,Dv(L_hX))L_h,\qquad
 F_h=h^{-n}(\det L_h)^2\widetilde f(L_hX),
\]
we have
$a_h(0)=0$, $c\le F_h(0)\le C$, and
\[
 [F_h]_{C^\alpha}\le C\|L_h\|^\alpha,\qquad
 [a_h]_{C^\beta}\le Ch^{-1}\|L_h\|^2
 \big(\|P_h\|\|L_h\|^\beta+\|P_h\|^2\|L_h\|\big)=o(1).
\]
The last assertion follows from \eqref{eq:22}, provided $\varepsilon>0$ is sufficiently small. We normalise by $F_h(0)^{1/n}$ and fix a small flatness
$E$ for Lemma~\ref{lem:14}. At a sufficiently small common height,
the coefficient seminorms are at most $\delta E$, and the quadratic
compactness yields a flatness of at most $E$. This lemma gives a uniform
rescaled Hessian estimate. At this fixed height, the coordinate matrices
and their inverses are uniformly bounded. The pullback, together with a finite covering argument,
yields the estimate on a fixed interior set overlapping the collar.
These estimates control all  pairs of points that are sufficiently close. The uniform
Hessian bound and the initial $C^1$ bound together complete the proof.
\end{proof}

\raggedbottom

\end{document}